\documentclass[11pt]{article}
\usepackage[english]{babel}
\usepackage[hmarginratio=3:2]{geometry}
\usepackage[titletoc]{appendix}
\usepackage{authblk}
\usepackage{graphicx}
\usepackage{amsmath}
\usepackage{amsthm}
\usepackage{amssymb}
\usepackage[mathscr]{eucal}
\usepackage{textcomp}
\usepackage[latin1]{inputenc}
\usepackage{latexsym}
\usepackage{hyperref}
\usepackage[dvipsnames]{xcolor}
\usepackage{framed}
\usepackage{tikz-cd}
\usepackage{ifsym}
\usepackage{enumerate}
\usepackage{sectsty}
\usepackage{chngcntr}

\numberwithin{equation}{section}

\numberwithin{subsection}{section}

\allowdisplaybreaks[1]

\newtheorem{theorem}{Theorem}[section]
\newtheorem{lemma}[theorem]{Lemma}
\newtheorem{corollary}[theorem]{Corollary}
\newtheorem{proposition}[theorem]{Proposition}
\newtheorem{observation}[theorem]{Observation}
\newtheorem*{remark}{Remark}
\newtheorem*{definition}{Definition}

\newenvironment{enumerate1}
{\begin{enumerate}[\upshape (1)]}
{\end{enumerate}}

 \newcommand\cB{\mathcal{B}}
\newcommand\cC{\mathcal{C}} \newcommand\cD{\mathcal{D}}
\newcommand\cE{\mathcal{E}} \newcommand\cF{\mathcal{F}}
\newcommand\cG{\mathcal{G}} 
\newcommand\cI{\mathcal{I}} 
 \newcommand\cL{\mathcal{L}}
 \newcommand\cN{\mathcal{N}}
\newcommand\cO{\mathcal{O}} 
 
\newcommand\cS{\mathcal{S}} 
\newcommand\cU{\mathcal{U}} \newcommand\cV{\mathcal{V}}
 \newcommand\cX{\mathcal{X}}
\newcommand\cY{\mathcal{Y}}

\newcommand\GG{\mathbb{G}}

 \newcommand\NN{\mathbb{N}}
 \newcommand\PP{\mathbb{P}}
\newcommand\QQ{\mathbb{Q}}

 \newcommand\ZZ{\mathbb{Z}}

 \newcommand\rH{\mathrm{H}}
 
\newcommand\rK{\mathrm{K}}

\newcommand\rQ{\mathrm{Q}} \newcommand\rR{\mathrm{R}}
\newcommand\rS{\mathrm{S}}

\newcommand\rmm{\mathrm{m}}

\newcommand\bfg{\mathbf{g}}

 \newcommand\frm{\mathfrak{m}}

\newcommand\arr{\ifinner\to\else\longrightarrow\fi}

\newcommand\arrto{\ifinner\mapsto\else\longmapsto\fi}

\newcommand{\xarr}{\xrightarrow}

\newcommand{\eqdef}{\mathrel{\smash{\overset{\mathrm{\scriptscriptstyle def}} =}}}

\newcommand\into{\hookrightarrow}

\renewcommand{\th}{^\text{th}}

\def\displaytimes_#1{\mathrel{\mathop{\times}\limits_{#1}}}

\def\displayotimes_#1{\mathrel{\mathop{\bigotimes}\limits_{#1}}}

\renewcommand\hom{\operatorname{Hom}}

\newcommand\aut{\operatorname{Aut}}

\newcommand\spec{\operatorname{Spec}}

\newcommand\rk{\operatorname{rk}}

\newcommand\pr{\operatorname{pr}}

\renewcommand\projlim{\varprojlim}

\newcommand{\underhom}{\mathop{\underline{\mathrm{Hom}}}\nolimits}

\newcommand{\underaut}{\mathop{\underline{\mathrm{Aut}}}\nolimits}

\newlength{\ignora}

\renewcommand\projlim{\varprojlim}

\newcommand{\ind}{\operatorname{Ind}}

\newcommand{\mmu}{\boldsymbol{\mu}}

\newcommand{\gm}{\GG_{\rmm}}

\newcommand{\GL}{\mathrm{GL}}

\newcommand{\abs}[1]{\left|#1\right|}

\newcommand{\step}[1]{\smallskip\textit{#1.}\;}

\renewcommand{\epsilon}{\varepsilon}

\renewcommand{\k}[1][*]{\operatorname{K}'_{#1}}

\newcommand{\kz}{\operatorname{K}_{0}}

\newcommand{\R}{\mathrm{R}}
\newcommand{\Rt}{\smash{\widetilde{\mathrm{R}}}}

\newcommand{\Q}[1]{\QQ(\zeta_{#1})}

\newcommand{\kmu}[2]{\rK'_{*}(#2)_{(\mmu_{#1})}}

\newcommand{\kg}[2][*]{\rK'_{#1}(#2)_{\bf g}}
\newcommand{\ka}[2][*]{\rK'_{#1}(#2)_{\bf a}}
\newcommand{\kt}[2][*]{\rK'_{#1}(#2)_{\bf t}}
\newcommand{\knt}[2][*]{\rK'_{#1}(#2)_{\bf nt}}

\newcommand{\lc}{\overline{\mathcal{C}}}

\newcommand{\red}{_{\rm red}}

\newcommand{\ci}[1][\relax]{\cI_{\mmu_{{#1}}}}
\newcommand{\cd}[1][\relax]{\cD_{\mmu_{{#1}}}}

\newcommand{\rh}[2]{\rH_{#1}(#2)}
\newcommand{\rhin}[2]{\rH^{\rm in}_{#1}(#2)}
\newcommand{\rhmu}[2]{\rH_{\mmu_{#1}}(#2)}
\newcommand{\rhmuin}[2]{\rH^{\rm in}_{\mmu_{#1}}(#2)}

\newcommand{\cha}{\operatorname{char}}

\newcommand{\autm}[1][\infty]{\aut{\mmu_{#1}}}

\newcommand{\coh}{\operatorname{\mathfrak{Coh}}}

\newcommand\nome{testing}
\newcommand\call[1]{\label{#1}\renewcommand\nome{#1}}
\newcommand\itemref[1]{\item\label{\nome;#1}}
\newcommand\refall[2]{\ref{#1}(\ref{#1;#2})}
\newcommand\refpart[2]{(\ref{#1;#2})}

\title{On a Riemann-Roch formula for stacks with finite cyclotomic inertia}
\author[1]{Francesco Sala}
\affil[1]{\textit{Scuola Normale Superiore, Piazza dei cavalieri 7, 56126, Pisa, Italy}}
\date{}

\begin{document}
\maketitle

\begin{abstract}
In \cite{To}, B. Toen defined a Riemann-Roch map from the rational algebraic K-theory of a tame Deligne-Mumford quotient stack to the étale K-theory of its inertia. He proved that this map is covariant with respect to proper maps. Moreover in \cite{vezzosi-vistoli02} G. Vezzosi and A. Vistoli proved a decomposition theorem for the equivariant K-theory of a noetherian scheme. In this paper we generalize Toen's construction and give a geometric definition of the Vezzosi-Vistoli decomposition, interpreting the pieces as corresponding to the components of the \emph{cyclotomic inertia}. When the map from the cyclotomic inertia to the stack is finite, we can define a Riemann-Roch map in Toen's style. We prove that this map is an isomorphism and it is covariant with respect to proper relatively tame maps; moreover in some favourable circumstances we explicitly compute its inverse map, and show that we can recover Toen's one when the stack is tame Deligne-Mumford.
\end{abstract}

\bigskip
\bigskip
In \cite{To}, B. Toen defined a \emph{Riemann-Roch map} from the rational algebraic K-theory of a tame Deligne-Mumford quotient stack to the étale K-theory of its inertia tensored with the field $\QQ(\zeta_\infty)$ containing all roots of unity. Then he proved that his map is covariant with respect to proper-push-forwards.\\

In this paper we extend the work of Toen to more general settings, treating for example the case where the stack is tame in the sense of \cite{dan-olsson-vistoli} (not necessarily Deligne-Mumford). We make use of a variant of the inertia stack, that is the \emph{cyclotomic inertia}; in characteristic $0$ this was defined in \cite{AGV}, and is a twisted form of the classical inertia stack. Here we give a general definition of the cyclotomic inertia, valid in every characteristic, following the work of L. Schadeck (\cite{ScT}). 

We identify a nice class of algebraic stacks with finite stabilizers, that is those such that the stack $\cX:=[X/G]$ and its cyclotomic inertia $\cI_{\mu}\cX$ are related by a finite map $\rho\colon \cI _{\mu}\cX\rightarrow \cX$; this comprehend tame stacks and DM stacks over an equal characteristic base ring.

In this case the proper push-forward $\rho_*$ makes it possible to define a covariant \emph{Lefschetz-Riemann-Roch map}, whose properties we will explore.\\

The codomain of this map will be a twisted form of the \emph{geometric K-theory}; the latter is a suitable localization of the rational algebraic K-theory of the stack, first defined in \cite{Schadeck}, that is inspired by the work of Vistoli-Vezzosi summarized below.
 
By \emph{twisting}, we mean that we will consider a suitable tensor $\kg{\cI_{\mmu_r}\cX}\otimes \QQ(\zeta_r)$, depending on the connected components of $\cI_{\mmu}\cX$. This formulation differs slightly from Toen's one, since in this case the "coefficients" $\QQ(\zeta_r)$ are not invariant by push-forward. Anyhow, we will prove a comparison theorem with respect to Toen's result.

In doing so, we prove an existence theorem for rigidifications and Chow envelopes for tame stacks, that may be of independent interest.\\
 
In their paper \cite{vezzosi-vistoli02}, G. Vezzosi and A. Vistoli proved a decomposition theorem for the rational equivariant K-theory of an algebraic space $X$ by an affine algebraic group $G$. The pieces of this decomposition are indexed by the conjugacy classes of subgroups of $G$ which are isomorphic to the groups $\mmu_n$ of multiplicative type, called \emph{dual cyclic subgroups}.

Assuming a technical hypothesis (that the action is \emph{sufficiently rational}, see the introduction of \cite{vezzosi-vistoli02}), they proved the following theorem: defining $\tilde R(\sigma)$ to be $\QQ(\zeta_n)$ when $\sigma$ has order $n$ and $w_G(\sigma):=N_G(\sigma)/C_G(\sigma)$, where $N_G(\sigma),C_G(\sigma)$ are the normalizer and the centralizer of $\sigma$ in $G$: 
\[
K(X,G)\simeq \underset{\sigma\in \cC(G)}{\prod} (K(X^\sigma, C_G(\sigma))_{geom}\otimes \tilde R(\sigma))^{w_G(\sigma)}
\]
In the formula, $K(-)_{geom}$ is the localization of the rational equivariant K-theory by the multiplicative subset of elements of non-zero rank. \\

In this paper we interpret the dual cyclic subgroups of Vezzosi-Vistoli as indexing the connected components of the cyclotomic inertia; in particular we will be able to describe the functorial behaviour of the Vezzosi-Vistoli decomposition in a wide range of cases. Moreover when $\cX$ is a global quotient by a finite group scheme or tame or Deligne-Mumford we will define an inverse to our Riemann-Roch map, and show that in the tame DM case over a field it agrees with Toen's map. Finally we show that this gives a (non canonical) isomorphism from the rational K-theory of a quotient stack to the Chow group of its cyclotomic inertia.

\subsection*{Summary of the paper}
In the first section we define the \emph{cyclotomic inertia} $\cI\cX=\sqcup_r \cI_{\mmu_r}\cX$ of an algebraic stack $\cX$ with finite stabilizers, and  in the case where $\cX=[X/\GL_n]$ we give an explicit description of it in terms of fixed loci of subgroups of $\GL_n$ on $X$. Moreover, we prove that the construction is functorial (that is, a map of stacks induces a map between the corresponding cyclotomic inertias).\\

In the third section we recall the main results of \cite{Schadeck}; the most important one for the current paper will be an intrinsic definition of the \emph{geometric K-theory} of a stack $\kg{\cX}$, defined as a localization of the algebraic K-theory $\k(\cX)$; meanwhile, we see that there is a decomposition $\k(\cX)=\oplus_r \k(\cX)_{(r)}$, the geometric part corresponding to $r=1$. Identifying $\k(\cX)$ with an equivariant K-theory ring $K(X,G)$ (for $\cX=[X/G]$) this construction coincides with the Vezzosi-Vistoli one in \cite{vezzosi-vistoli02}; our result ensures that this construction does not depend on the presentation  $\cX=[X/G]$.\\

In the fourth section we introduce the \emph{tautological part} $\kt{\cI\cX}$ of the K-theory of the cyclotomic inertia stack, and we prove that it is isomorphic to the algebraic K-theory of the original stack $\k(\cX)$. In particular this will allow us to give a Riemann-Roch map landing in the tautological K-theory of the cyclotomic inertia.\\

In the fifth section we define the twisting map $\alpha\colon \bigoplus_r \kg{\cI\cX}\otimes\QQ(\zeta_r) \arr \kt{\cI\cX}$ and prove that it is covariant. This will allow us to define a \emph{Lefschetz-Riemann-Roch} map with values in $\bigoplus_r \kg{\cI\cX}\otimes\QQ(\zeta_r)$.\\

In the sixth section we prove that the twist map is an isomorphism; this allows us to put togheter all the ingredients of the previous sections and define a Lefschetz-Riemann-Roch map. In the cases where $\cX$ is a global quotient by a finite group, is tame or Deligne-Mumford we give an explicit formula for this map, similar to the one given by Toen.

Finally, when $\cX$ is tame Deligne-Mumford we prove that our formula can be used to recover Toen's one (up to identifying geometric and étale K-theory of a stack). The main result is a comparison of $\bigoplus_r \kg{\cI\cX} \otimes\QQ(\zeta_r)$ with $\kg{\cI\cX}\otimes \QQ(\zeta_\infty)$ (where in the lhs the coefficients $\QQ(\zeta_r)$ are not constant after push-forwards, while in the rhs they are constant).\\

In the last section, finally, we prove that there is a Riemann-Roch map with values in $\kg{\cI\cX}$ and rational coefficients. This map is noncanonical, but it allows nonetheless to define a Chern character with values in the rational Chow ring of the cyclotomic inertia.
\subsection*{Acknowledgements}

This paper is part of the author's Phd thesis. The project was proposed to me by my advisor A. Vistoli, and I would like to thank him sincerely for sharing with me his insights and for patiently accompanying me through the various stages of this work. 

\subsection*{Notation}

In what follows we will assume that $\cX$ is a separated algebraic stack with finite stabilizers, and that it can be written as a quotient $\cX=[X/G]$ where $X$ is a noetherian scheme over an affine connected base scheme $A=Spec(R)$, with $R$ excellent, and $G=GL_{n,R}$. 
The last condition is actually equivalent to the apparently weaker one that $\cX=[X/G]$, where $G$ is a flat $R-$linear algebraic group: indeed by choosing an embedding $G \subseteq \GL_{n} \eqdef \GL_{n,R}$, and replacing $X$ with $X\times^{G}\GL_{n}$, we may suppose that $G = \GL_{n}$. We will always write $\cX=[X/G]$, where $G$ is a product of general linear groups on $R$.

As in \cite{vezzosi-vistoli02}, we call a subgroup scheme $\sigma \subseteq G$ \emph{dual cyclic} when it is isomorphic to $\mmu_{r, R}$ for some $r \geq 1$.

For each $r \geq 1$, the datum of a monomorphism from a dual cyclic subgroup $\mmu_{r,R}\arr\GL_{n,R}$ is equivalent to a locally free splitting $R^n=\bigoplus_{i=1}^r S_i$; the equivalence is obtained considering the decomposition of $R^n$ into eigenspaces for the action of $\mmu_{r,R}$.

\begin{remark}
Consider a map $\mmu_{r,R}\arr\GL_{n,R}$ and the corresponding decomposition $R^n=\bigoplus_{i=1}^r S_i$. The $S_i$'s are not necessarily free, and in those cases we say that the cyclic subgroup is not constant; otherwise we call it constant.
\end{remark}

Let $\overline{\cC}_r(G)$ be the set of \emph{conjugacy classes} of monomorphisms of dual cyclic subgroups of $G$ of order $r$. More precisely, the functor which assigns to an $R-$scheme $S$ the set of splittings $R^n=\bigoplus_{i=1}^r S_i$ is represented by an open subscheme $U$ of a product of grassmanians (we shall elaborate more on this in the following sections), and it has a natural $G$-action (the "change of coordinates" action). The connected components of the quotient $[U/G]$ are indexed by the set of partitions $\nu= (\nu_1,\dots, \nu_r)$ of $n$ in $r$ parts (where each $\nu_i$ is the dimension of the $i$th eigenspace of $\mmu_{r,R}$); for any $\sigma$ we will call the corresponding $\nu$ its \emph{type}, and  $\overline{\cC}_r(G)$ will indicate the set of types. The moduli space of each $\nu$-component is isomorphic to $A$; indeed locally on $A$ every dual cyclic subgroup is constant, and conjugacy classes of constant dual subgroup schemes are determined by their \emph{type}.

If $G=\GL_{n_1}\times\cdots\times\GL_{n_k}$ is a product of general linear groups then the same discussion applies, with the only difference that the types are indexed by $r-$partitions of each set $\{1,\dots,n_i\}$.\\

Let $\overline{\cC}_r(G)$ be the set of \emph{types} of dual cyclic subgroups of $G$, or equivalently the set of constant dual cyclic subgroups of $G$. From now on, when we pick a dual cyclic subgroup $\sigma\subset G$ of any given type, we will mean a \emph{constant} representative of the type.

The group $Aut(\mmu_r)$ has a natural action on $\overline{\cC}_r(G)$, and we let $\cC_r(G)$ be the set of orbits for this action.\\

Let $\k(\cX)$ denote the $K-$theory of coherent sheaves on $\cX$, while let  $\kz$ denote the naive ring of locally free sheaves on $\cX$.

Finally, for a group scheme $G$ of multiplicative type over $A$, let $\R G$ be the ring of characters of $G$. If $G=\GL_{n_1}\times\dots\times\GL_{n_k}$, let $\R G$ be the Grothendieck ring of highest-weight algebraic representations of $G$ (which is isomorphic to $\QQ[\hat{T}]^W$, where $\hat T$ denotes the characters of a maximal torus and $W$ is the Weyl group of $G$).\\

There is a special class of stacks and morphisms for which the properties of the decomposition discussed above are stronger: this is the class of \emph{tame} algebraic stacks: following \cite{dan-olsson-vistoli}, such a stack is called tame if the stabilizer of any point is linearly reductive. 

In general, a morphism of stacks $f\colon\cX\rightarrow \cY$ is called relatively tame if the relative inertia $\cI_{\cY}\cX$ has linearly reductive geometric fibers (see \cite{dan-olsson-vistoli}). For such a morphism, the push-forward $f_*$ descends to K-theory.

\section{The cyclotomic inertia of a separated stack}

In this section we define an explore the notion of \emph{cyclotomic inertia}. This notion was first defined in \cite{AGV} in order to define Gromov-Witten invariants for Deligne-Mumford stacks.

In the case wehere $S$ is the spectrum of a field most of these properties have already been defined by A. Vistoli and L. Schadeck, and can be found in Schadeck's Phd thesis (\cite{ScT}).\\

Let $\cX$ be as above, $r$ be a positive integer. The \emph{$r\th$ cyclotomic inertia} is a stack $\ci[r]\cX$ fibered in groupoids over the category of schemes over $R$, defined as follows.

An object of $\ci[r]\cX(S)$, where $S$ is a scheme, is a pair $(\xi, a)$, where $\xi$ is an object of $\cX(S)$ and $a\colon \mmu_{r, S} \arr \underaut_{S}\xi$ is a monomorphism of group schemes. An arrow $f\colon (\xi, a) \arr (\eta, b)$ from $(\xi, a) \in \ci[r]\cX(S)$ to $(\eta, b) \in \ci[r]\cX(T)$ consists of an arrow $f\colon \xi \arr \eta$ in $\cX$, such that the diagram
   \[
   \begin{tikzcd}
   \mmu_{r, S} \ar[r, equal]\ar[d, "a"] & \phi^{*}\mmu_{r, T}\ar[d, "\phi^{*}b"]\\
   \underaut_{S}\xi \ar[r]& \phi^{*}\underaut_{T}\eta
   \end{tikzcd}
   \]
commutes. Here $\phi\colon S \arr T$ is the image of $f$ in the category of schemes, and the bottom row is induced by $f$.

There is an obvious morphism of fibered categories $\ci[r]\cX \arr \cX$ that sends $(\xi, a)$ into $\xi$.

For the main properties of this construction, we refer to \cite{Schadeck}. Here we will just mention the relevant definitions and facts.\\

Suppose that $\Delta$ is a finite diagonalizable group scheme over $R$, and $G$ is an affine group scheme of finite type over $R$. Consider the contravariant functor $\rh{\Delta}{G}$ from schemes over $R$ to sets, sending a scheme $S$ into the set of homomorphisms of group schemes $\Delta_{S} \arr G_{S}$. We will think of $\rh{\Delta}G$ as a Zariski sheaf. 
We denote by $\rhin{\Delta}G \subseteq \rh{\Delta}G$ the subsheaf consisting of the set of monomorphisms of group schemes $\Delta_{S} \arr G_{S}$.

If $\phi\colon \Delta \arr \Delta'$ is a homomorphism of diagonalizable group schemes, composing with $\phi$ gives a natural transformation $\rh{\Delta'}G \arr \rh{\Delta}G$. Call $\rQ(\Delta)$ the set of quotients of $\Delta$. For each $\Delta' \in \rQ(\Delta)$, consider the composite $\rhin{\Delta'}G \subseteq \rh{\Delta'}G \arr \rh{\Delta}G$, which is immediately seen to be a monomorphism. This induces a $G$-equivariant morphism of Zariski sheaves
   \[
   \coprod_{\Delta' \in \rQ(\Delta)} \rhin{\Delta'}G \arr  \rh{\Delta}G\,.
   \]

\begin{proposition}\call{prop:cyclotomic-in-G}\hfil
\begin{enumerate1}

\itemref{1} The functors $\rh{\Delta}{G}$ and\/ $\rhin{\Delta}G$ are represented by a quasi-projective scheme over $R$.

\itemref{2} If $\Delta' \in \rQ(\Delta)$, then $\rhin{\Delta'}G$ is open and closed in $\rh{\Delta}G$.

\itemref{3} The morphism $\coprod_{\Delta' \in \rQ(\Delta)} \rhin{\Delta'}G \arr  \rh{\Delta}G$ is an isomorphism.

\itemref{4} If $G$ is finite and linearly reductive, then $\rh{\Delta}{G}$ and $\rhin{\Delta}G$ are finite over $\spec R$. If $G$ is étale and every localization $R_p$ is equicharacteristic,  $\rh{\Delta}{G}$ and $\rhin{\Delta}G$ are finite over $\spec R$.


\end{enumerate1}

\end{proposition}
The proof can be found in the appendix, Proposition \ref{prop:cyclotomic-in-G}.\\

There is an action of $G$ on $\rh{\Delta}G$, by conjugation, leaving $\rhin{\Delta}G$ invariant. Consider the closed subscheme $X_{\mmu_{r}} \subseteq X \times \rhmuin{r}G$ defined as follows. Let $(x, \phi)$ be a point of $\bigl(X \times \rhmuin{r}G\bigr)(S)$; that is, $x$ is a morphism of schemes $S \arr X$, and $\phi\colon \mmu_{r, S} \arr G_{S}$ is a monomorphism of group schemes. Then $\phi$ induces an action of $\mmu_{r, S}$ on the $S$-scheme $S \times X$; we say that $(x, \phi)$ is in $X_{\mmu_{r}}(S)$ if the section $S \arr S \times X$ defined by $x$ is fixed under the action of $\mmu_{r, S}$. 

The subscheme $X_{\mmu_{r}} \subseteq X \times \rhmuin{r}G$ is  $G$-invariant; the projection $X_{\mmu_{r}} \arr X$ is $G$-equivariant, and defines a morphism of algebraic stacks $[X_{\mmu_{r}}/G] \arr [X/ G] = \cX$.

\begin{proposition}\label{prop:cyclotomic-inertia-quotient}
The quotient stack $[X_{\mmu_{r}}/G]$ is naturally equivalent to $\ci[r]\cX$. 
\end{proposition}

\begin{proof}
We can immediately exhibit an equivalence. Namely, let $(\xi,a)$ be a $S-$point of $\ci[r]\cX$: then $\xi$ correspond to an equivariant morphism $g\colon P\arr X$, where $P$ is a principal $G-$bundle over $S$; moreover $a$ gives a morphism $\mmu_{r,S}\arr G_S$ such that $g$ factors through $P/\mmu_{r,S}\arr X$.

By definition these data give a unique equivariant morphism $P\arr  X \times \rhmuin{r}G$, and it is immediate to see that it factors through $X_{\mmu_r}$. In particular they provide a well-defined $S-$point of $[X_{\mmu_{r}}/G]$, and we immediately see that this map can be extended to the morphisms $(\xi,a)\arr (\eta, b)$, giving a stack map.

Conversely, given a principal $G-$bundle $P\arr S$ and an equivariant map $P\arr X_{\mmu_{r}}$ we can compose it with the inclusion $X_{\mmu_{r}} \subseteq X \times \rhmuin{r}G$, giving two maps $\xi\colon P\arr X$ and $P\arr \rhmuin{r}G$. From the second one, quotienting by $G$, we get a map $\mmu_{r,s}\arr G_S$; moreover, $P\arr X$ being $\mmu_{r,P}$-invariant, we have a map $P/\mmu_{r,S}\arr X$ and $\mmu_{r,S}\arr G_S$ factors through $\underaut_{S}\xi$. This gives a morphism $[X_{\mmu_{r}}/G]\arr\ci[r]\cX$, and we immediately see that it is inverse to the above map.
\end{proof}

From this we can deduce the following important property.

\begin{proposition}\label{prop:inertia->stack}
Let $\cX$ be an algebraic stack over $R$. Then the morphism $\ci[r]\cX\arr \cX$ is representable. Moreover
\begin{enumerate1}

\item If $\cX$ is tame, then $\ci[r]\cX$ is also a tame algebraic stack; furthermore, the morphism $\ci[r]\cX\arr \cX$ is finite.
\item If $\cX$ is Deligne-Mumford and $R$ is equicharacteristic, then  $\ci[r]\cX$ is also Deligne-Mumford; furthermore, the morphism $\ci[r]\cX\arr \cX$ is finite.
\end{enumerate1}
\end{proposition}

\begin{proof}
Let us show that, in both cases, the morphism $\ci[r]\cX\arr \cX$ is representable and finite. Let $M$ be the moduli space of $\cX$: formation of $\ci[r]\cX$ commutes with base change on $M$, so the question is fppf local on $M$. Locally on $M$ the stack $\cX$ is of the form $[X/\Gamma]$, where $X \arr M$ is a finite morphism and $\Gamma \arr M$ is a linearly reductive finite group scheme (in the first case, \cite{dan-olsson-vistoli}) or an étale group scheme (in the second case, \cite{LM}, Chapter 6); by further refining in the fppf topology, we can assume that $\Gamma$ is obtained by pullback from a finite linearly reductive group scheme $G$ over $A$. From the construction above we have a factorization
   \[
   \ci[r]\cX \subseteq [(X \times \rhmuin{r}{G})/G] \arr [X/G] = \cX
   \]
where the first homomorphism is a closed embedding, and the second is finite because of Proposition~\refall{prop:cyclotomic-in-G}{4}. 
\end{proof}

\begin{remark}
The proofs in the above paragraph show, in particular, the following: if $\cX=[X/G]$ with $G=\GL_{n,R}$ then $\ci[r]\cX$ is the disjoint union
\[
\ci[r]\cX=\coprod_{\substack {\nu\in \lc_r(G)}}[X^\sigma/C_{G}(\sigma)].
\]
where $\sigma$ is any dual cyclic subgroup of  type $\nu$.

Indeed, for any $S$ the set $ \rhmuin{r}{G}(S)$ is parametrized by the dual cyclic subgroups of $\GL_n(S)$. For each such $\sigma$, $G$ acts on the set $(X^\sigma \times H_\mu{G})(S)$ with stabilizer $C_G(\sigma)$ (where $H_\mu$ is the subscheme of $\rhmuin{r}{G}$ parametrizing subgroups of type $\mu$). The result follows. 
\end{remark}

\subsection*{Functoriality of cyclotomic inertia stacks}

If $f\colon \cY \arr \cX$ is a representable morphism of stacks, there is an obvious induced morphism $\ci[r]f\colon \ci[r]\cY \arr \ci[r]{\cX}$: an object $(\eta, b)\in \ci[r]\cY$ is sent into $(f\eta,  f_{\eta}\circ b) \in \ci[r]\cX$, where $f_{\eta} \colon \aut_{S}\eta \arr \aut_{S}(f\eta)$ is the homomorphism induced by $f$. The point here is that $f_{\eta}$ is a monomorphism, because $f$ is representable. If $f$ is not representable then $f_{\eta}$ is not a monomorphism, so $(f\eta, f_{\eta}\circ b)$ is not an object of $\ci[r]\cX(S)$.

However, while the single $r\th$ inertia stacks are not functorial for nonrepresentable morphisms, their disjoint union is functorial. Let us define the \emph{cyclotomic inertia} $\ci\cX$ of $\cX$ as the disjoint union of the $\ci[r]\cX$. Since there a finite number of $r$ such that $\ci[r]\cX \neq \emptyset$, the morphism $\ci\cX \arr \cX$ is finite.

The stack $\ci\cX$ has an alternate description that makes its functoriality properties evident. Set $\mmu_{\infty} \eqdef \projlim_{n}\mmu_{n}$, where the index set is the set of positive integers ordered by divisibility, and homomorphisms $\mmu_{n} \arr \mmu_{m}$ are defined as $z \arrto z^{n/m}$. Then $\mmu_{\infty}$ is a affine group scheme over $k$, which is not of finite type.

Denote by $\ci[\infty]\cX$ the stack that is defined as follows. An object $(\xi, \phi)$ of $\ci[\infty]\cX$ over a scheme $S$ consists of an object $\xi$ of $\cX(S)$ and a homomorphism of group schemes $\phi\colon \mmu_{\infty, S} \arr \underaut_{S}\xi$. The arrows are defined in the obvious way. Forgetting $\phi$ gives a morphism $\ci[\infty]\cX \arr \cX$.

There is an morphism $\ci[r]\cX \arr \ci[\infty]\cX$ defined as follows. If $(\xi, a)$ is an object of $\ci[r]\cX(S)$, we obtain an object $(\xi, \phi)$ of $\ci[\infty]\cX(S)$, where $\phi$ is the composite of $a\colon \mmu_{r, S} \into \underaut_{S}\xi$ with the projection $\mmu_{\infty, S} \arr \mmu_{r,S}$. The definition for arrows is clear.

These give a morphism $\ci\cX \arr \ci[\infty]\cX$.

\begin{proposition}
The morphism $\ci\cX \arr \ci[\infty]\cX$ defined above is an equivalence of stacks over $\cX$.
\end{proposition}

\begin{proof}
Choose a positive integer $N$ that is divisible by the orders of the automorphism group schemes of all the geometric point of $\cX$. Then every homomorphism $\mmu_{\infty, S} \arr \underaut_{S}\xi$ factors uniquely through $\mmu_{N, S}$. Hence we can restate the definition of $\ci[\infty]\cX$ by taking objects $(\xi, \phi)$, where $\xi$ is an object of $\cX(S)$ and $\phi\colon \mmu_{N, S} \arr \underaut_{S}\xi$ is a homomorphism (this will actually give a stack that is strictly isomorphic to $\ci[\infty]\cX$).

This gives a description of the stack $\ci[\infty]\cX$ as a quotient stack, similar to the one we have given for $\ci[r]\cX$. Write $Y \subseteq X \times \rh{\mmu_N}G$ for the closed locus defined by the subfunctor of pair $(\xi, \phi)$ such that $\xi$ is fixed under the action of $\mmu_{N}$ defined by $\phi$; then $\ci[\infty]\cX = [Y/G]$.

For each positive integer $r$ that divides $N$, consider $\mmu_{r}$ as a quotient of $\mmu_{N}$; there is an embedding $\rhmuin{r}G \subseteq\rhmu{N}G$. Furthermore, the intersection of $Y$ with $X \times \rhmuin{r}G \subseteq X \times \rhmu{N}G$ is exactly the subscheme $X_{\mmu_{r}} \subseteq X \times \rhmuin{r}G$ that intervenes in the proof of Proposition~\ref{prop:cyclotomic-inertia-quotient}, so that $[X_{\mmu_{r}}/G] = \ci[r]\cX$, and the morphism $[X_{\mmu_{r}}/G] \subseteq [Y/ G]$ is isomorphic to the morphism $\ci[r]\cX \arr \ci[\infty]\cX$ defined above. The the result follows from Proposition \refall{prop:cyclotomic-in-G}{3}.
\end{proof}

From now on we will identify $\ci\cX$ with $\ci[\infty]\cX$, via the equivalence above.

This description of $\ci\cX$ makes its functorial property evident. If $f\colon \cY \arr \cX$ is a morphism of stacks and $(\eta, b)$ is an object of $\ci\cX(S) =\ci[\infty]\cY(S)$, then $(f\eta, b \circ f_{\eta})$ is in fact an object of $\ci[\infty]\cX(S)$. This function from the objects of $\ci\cY$ to the objects of $\ci\cX$ extends easily to a functor $\ci f\colon \ci\cY \arr \ci\cX$, such that the diagram
   \[
   \begin{tikzcd}
   \ci\cY \ar[r, "\ci f"]\ar[d] & \ci\cX\ar[d]\\
   \cY \ar[r, "f"] & \cX
   \end{tikzcd}
   \]
is strictly commutative. This gives a strict $2$-functor from stacks to stacks, sending $\cX$ into $\ci\cX$.

Yet another description of $\ci[r]\cX$ and $\ci\cX$ is as follows. Suppose that $(\xi, a)$ is an object of $\ci\cX(S)$; here $\xi$ is an object of $\cX(S)$, and $a\colon \mmu_{\infty, S} \arr \underaut_{S}\xi$ is a homomorphism of group schemes. This gives an action of $\mmu_{\infty}$ on $\xi$, which, in turn, gives a morphism of stacks $\phi\colon \cB_{S}\mmu_{\infty, S} \arr \cX$. Conversely, such a morphism of stacks gives an object $\xi$ of $\cX(S)$, defined as the image of the trivial torsor $\mmu_{\infty, S} \arr S$ through $\phi$, while the action of $\mmu_{\infty, S}$ on $\xi$ is defined as the homomorphism $\mmu_{\infty, S} = \underaut_{S}\mmu_{\infty, S} \arr \underaut_{S}\xi$. This gives an equivalence between $\ci\cX$ and the stack $\underhom_{k}(\cB_{k}\mmu_{\infty}, \cX)$, whose objects over $S$ are morphisms $\cB_{S}\mmu_{\infty, S} \arr \cX$. An arrow from $\phi\colon \cB_{S}\mmu_{\infty, S} \arr \cX$ to $\psi\colon \cB_{S}\mmu_{\infty, S} \arr \cX$ consists a morphism of schemes $f\colon S \arr T$, and a $2$-commutative diagram
   \[
   \begin{tikzcd}[row sep = 3pt]
   \cB_{S}\mmu_{\infty, S} \ar[rd, "\phi"]\ar[dd] &\\
   & \cX\\
   \cB_{S}\mmu_{\infty, T} \ar[ur, swap, "\psi"]
   \end{tikzcd}
   \]
in which the column is the homomorphism $\cB_{S}\mmu_{\infty, S} \arr \cB_{T}\mmu_{\infty, T}$ induced by $f\colon S \arr T$.

A similar description can be given for $\ci[r]\cX$. Consider a morphism $\phi\colon \cB_{S}\mmu_{r, S} \arr \cX$; call $\xi \in \cX(S)$ the image of the trivial torsor, and $a\colon \mmu_{r,S} \arr \underaut_{S}\xi$ the induced homomorphism. Then $\phi$ is representable if and only if $a$ is an embedding of group schemes. Hence $\ci[r]\cX$ can be described as follows: the objects of $\ci[r]\cX$ consist of representable morphisms $\cB_{S}\mmu_{r, S} \arr \cX$. The arrows can be described as in the case of $\ci\cX$.

For each $r > 0$ consider the automorphism group $\autm[r]$ of $\mmu_{r}$, which is the group of units in the ring $\ZZ/r\ZZ$. There is obvious strict right action of $\autm[r]$ on $\ci[r]\cX$, in which an automorphism $u$ of $\mmu_{r}$ sends an object $(\xi, a)$ of $\ci[r]\cX(S)$ into $(\xi, a \circ u)$. 

Consider now the profinite group $\autm = \projlim_{r}\autm[r]$; this acts on each $\ci[r]\cX$ through its quotient $\autm[r]$. This induces an action of $\autm$ on $\ci\cX = \coprod_{r}\ci[r]\cX$. This can also be described using the identification of $\ci\cX$ with $\ci[\infty]\cX$ by the same formula as above.

The stack-theoretic quotient $[\ci[r]\cX/\autm[r]]$ can be described as follows. Denote by $\cd[r]\cX$ the stack whose object over $S$ are triples $(\xi, \sigma, a)$, where $\xi$ is an object of $\cX(S)$, $\sigma \arr S$ is a finite group scheme, étale-locally isomorphic to $\mmu_{r, S}$, and $a\colon \sigma \arr \underaut_{S}\xi$ is an embedding of group schemes over $S$. The arrows are defined in the obvious way.

There is a natural morphism $\ci[r]\cX \arr \cd[r]\cX$ that associates with each object $(\xi, a)$ of $\ci[r]\cX(S)$ the object $(\xi, \mmu_{r,S}, a)$ of $\cd[r]\cX$. This is easily seen to make $\ci[r]\cX$ into a $\autm[r]$-torsor over $\cd[r]\cX$; the corresponding action of $\autm[r]$ on $\ci$ is exactly the one described above. This construction will not be used in what follows.

\section{Conventions}

Let us recall the notation. Let $A=Spec(R)$ be a noetherian excellent connected scheme; all the objects in the subsequent sections will be intended relative to $A$.

In particular let $\cX=[X/G]$ be a separated algebraic stack, with $G=\GL_{n_1}\times\dots\times\GL_{n_k}$ a finite product of general linear groups. Moreover, we will always assume that the cyclotomic inertia map $\ci\cX \arr\cX$ is finite: this is fundamental to compare the K-theory of $\ci\cX$ and that of $\cX$ via the push-forward map. 

It is easy to see that the cyclotomic inertia is not finite over $\cX$ whenever there is a stabilizer diagonalizable group scheme degenerating to a non-diagonalizable group over a closed subset: for example when $\cX=B_RG$ with $R$ a dvr such that $G$ is diagonalizable over the generic fiber and unipotent over the special fiber. We may consider $G=\ZZ\slash p\ZZ$ over $Spec(\ZZ_p)$ or $G=Spec(\mathbb F_p[x,y]_{(x)}\slash(y^p-1))$ as a group scheme over $Spec(\mathbb F_p[x]_{(x)})$ with group law $(y_1,y_2)\arr y_1+y_2+x\cdot y_1y_2$. In those cases the fiber of $\ci\cX \arr\cX$ over the closed point of $A$ is empty, while it is nonempty over the generic point of $A$: in particular $\ci\cX \arr\cX$ cannot be proper. \\

\section{The intrinsic decomposition of the K-theory\\of a quotient stack}

In this section we briefly recall the main notations and the results of \cite{Schadeck}, in particular the fundamental notion of \emph{geometric K-theory}. This will be fundamental in order to get a Riemann-Roch isomorphism.\\
 
Fix now $G=\GL_{n,R}$. For any $r \geq 1$ we have a decomposition $\R\mmu_{r} = \prod_{s|r}\Q{s}$. We denote by $\R \sigma\arr\Rt\sigma$ the projection onto the $\Q{r}$ factor.
Recall from \cite{vezzosi-vistoli02} that if $\sigma \subseteq G$ is a constant dual cyclic subgroup, $\frm_{\sigma} \subseteq \R G$ is the maximal ideal that is the kernel of the composite $\R G \arr \R\sigma \arr \Rt\sigma$; this only  depends on the \emph{type} of $\sigma$. 
If $M$ is an $\R G$-module, the \emph{$\sigma$-localization} $M_{\sigma}$ is the localization $M_{\frm_{\sigma}}$. If $\sigma \subseteq G$ is a constant dual cyclic subgroup, then $\k(X, G)_{\sigma} \neq 0$ if and only if $X^{\sigma} \neq \emptyset$.

We define the multiplicative system $\Sigma_{r}^{\cX} = \Sigma_{r} \subseteq \kz \cX$ as follows. An element $\alpha \in \kz\cX$ is in $\Sigma_{r}$ if for all representable morphisms $\phi\colon \cB_{S}\mmu_{r} \arr \cX$, where $S$ is an $R$-scheme, the projection of $\phi^{*}a \in \R\mmu_{r}$ in $\Rt\mmu_{r} = \Q{r}$ is non-zero.

In particular, $\Sigma_{1}$ admits the following description. If $\cX_{1}$, \dots,~$\cX_{m}$ are the connected components of $\cX$, then $\k\cX = \bigoplus_{i=1}^{m}\k\cX_{i}$. There is a rank map $\rk_{i}\colon \k\cX_{i} \arr \QQ$. Then $\alpha$ is in $\Sigma_{1}$ if and only if $\rk_{i}\alpha \neq 0$ for all $i$.

A morphism $f\colon \cY \arr \cX$ of stacks induces a pullback homomorphism $f^{*}\colon \kz\cX \arr \kz\cY$. If $f$ is representable, one immediately sees that $f^{*}$ carries $\Sigma_{r}^{\cX}$ into $\Sigma_{r}^{\cY}$. If $f$ is not representable this is certainly not true in general, but it is if $r = 1$. 

\begin{definition}
The $\mmu_{r}$-localization $\kmu r \cX$ of $\k\cX$ is the $\kz\cX$-module $\Sigma_{r}^{-1}\k(\cX)$.
\end{definition}

Set $\k(X,G)_{(r)} \eqdef \prod_{\sigma\in \cC_{r}(G)}\k(X, G)_{\sigma}$, so that $\k(X,G) = \prod_{r \geq 1} \k(X,G)_{(r)}$. The support of $\k(X, G)_{\sigma}$ in $\spec \R G$ is the maximal ideal $\frm_{\sigma} = \ker(\R G \arr \Rt\sigma)$.\\

Let $\sigma_{1}$, \dots,~$\sigma_{m}$ be (constant) representatives of dual cyclic subgroups $\sigma \in \cC(G)$  with $\abs{\sigma} = r$ such that $\k(X, G)_{\sigma} \neq 0$. Then
\[
   \k(X,G)_{(r)}=\prod_{\substack {\sigma\in \cC(G)\\ \abs{\sigma} = r}}\k(X, G)_{\sigma} =
   \prod_{i=1}^{m}\k(X, G)_{\sigma_{i}}.
\]

Suppose that $\phi\colon\cB_{S}\sigma \arr \cX$ is a representable morphism, where $\sigma$ is a constant dual cyclic group over $S$. We may assume that there is an embedding $\sigma \subseteq G_{S}$, and a rational point $p \in X(S)$ which is fixed under the action of $\sigma$. Hence $X_{S}^{\sigma} \neq \emptyset$, so $\sigma$ is conjugate to some $\sigma_{i}$ (that is it belongs to the same type). We say in this case that $\phi$ is \emph{associated} to $\sigma_i$.

Let $\sigma$ be any of the $\sigma_i$'s; we define another multiplicative system $\Sigma_\sigma^\cX=\Sigma_\sigma\subseteq \kz \cX$ as follows: $\alpha\in\kz\cX$ lies in $\Sigma_\sigma$ if and only if for any $S/A$ and any $\phi\colon \cB_S\mmu_r\arr \cX$ associated to $\sigma$ the pullback $\phi^*\alpha\in \R\mmu_r$ has a nonzero projection to $\Rt\mmu_r$.

\begin{proposition}\label{prop:decomposition}
Let $r$ be any positive integer and $\sigma \in \cC(G)$ be any constant dual cyclic subgroup such that $\k(X, G)_{\sigma} \neq 0$.
\begin{enumerate1}
\item The projection $\k(X, G) \arr \prod_{\sigma \in \cC(G)}\k(X, G)_{\sigma}$ is an isomorphism.
\item The $\mmu_{r}$-localization $\k \cX \arr \kmu r \cX$ isomorphic to the projection $\k(X,G) \arr \k(X,G)_{(r)}$.
\item The localization $\Sigma_\sigma^{-1}\k(\cX)$ is isomorphic to $\k(X,G)_\sigma$.
\end{enumerate1}

\end{proposition}

As a corollary, we get the following.

\begin{theorem}
The projections $\k(\cX) \arr \kmu r {\cX}$ induce an isomorphism
   \[
   \k(\cX) \simeq \prod_{r \geq 1}\kmu r {\cX}\,.
   \]
\end{theorem}

Another consequence is this.

\begin{corollary}\label{vanishing}
We have $\kmu r \cX \neq 0$ if and only if there exists a geometric point of $\cX$ whose automorphism group scheme contains a dual cyclic subgroup of order $r$.

In particular when $\cX$ is tame, $\kmu r \cX = 0$ for all $r > 1$ if and only if $\cX$ is an algebraic space.
\end{corollary}

In the decomposition $\k\cX = \prod_{r}\kmu r \cX$ we split off the factor $\kmu 1 \cX$, which we call the \emph{geometric K-theory} of $\cX$, and denote by $\kg\cX$, and the factor $\prod_{r \geq 2}\kmu r \cX$, which we call the \emph{algebraic part} of the K-theory of $\cX$, and denote by $\ka\cX$. Thus we have the \emph{fundamental decomposition}
   \[
   \k\cX = \kg\cX \oplus \ka\cX\,.
   \]

The fundamental decomposition $\k\cX = \kg\cX \oplus \ka\cX$ gives both a projection $\k\cX\arr\kg\cX$ and an embedding $\kg\cX \arr \k\cX$. 

\begin{proposition}\label{prop:functoriality}
Let $f\colon \cY \arr \cX$ be a homomorphism of stacks.

\begin{enumerate1}

\item Assume that $f$ has finite Tor dimension. Then there exists a homomorphism
   \[
   f^{*}\colon \kg\cX \arr \kg\cY
   \]
such that the diagram
   \[
   \begin{tikzcd}
   \k\cX \ar[r, two heads]\ar[d, "f^{*}"]& \kg\cX\ar[d, "f^{*}"]\\
   \k\cY \ar[r, two heads] & \kg\cY
   \end{tikzcd}
   \]
commutes.

\item Suppose that $f$ is proper and relatively tame. Then there exists a homomorphism
   \[
   f_{*}\colon \kg\cY \arr \kg\cX
   \]
such that the diagram
   \[
   \begin{tikzcd}
   \kg\cY \ar[r, hook]\ar[d, "f_{*}"]& \k\cY\ar[d, "f_{*}"]\\
   \kg\cX \ar[r, hook] & \k\cX
   \end{tikzcd}
   \]
commutes.

The maps $\kg- \into \k(-)$ and\/ $\k(-) \twoheadrightarrow \kg-$  in the diagrams above are the embeddings and the  projections coming from the fundamental decomposition. 
\end{enumerate1}
\end{proposition}

The homomorphism $f^{*}$ and $f_{*}$ defined above are clearly unique; they make $\kg-$ into a contravariant functor for maps of finite Tor dimension, and a covariant functor for proper maps.

\begin{remark}
\label{remark:push-forward}
We can make the above proposition more precise.Suppose that we have a morphism of stacks $f\colon [Z/G\times H]\arr [X/G]$ induced by an equivariant map $Z\arr X$.

Let $\rho$ be a dual cyclic subgroup of $G\times H$ and $\sigma$ be its projection to $G$. 

\begin{enumerate1}
\item The pull-back $f^*$ is a map of $\R G$-modules, so it induces a map
\[
f^*:\k(X,G)_\sigma=S_\sigma^{-1}\k(X,G)\arr S_\sigma^{-1}\k(Z,G\times H)=\prod_{\substack{\rho \in \cC(G\times H)\\
   \rho \twoheadrightarrow \sigma}}\k(Z, G \times H)_{\rho}
\]
\item The push-forward $f_*\colon \k(Z,G\times H)\arr \k(X,G)$ is a morphism of $\R G$-modules. If $\eta\in \k(Z, G\times H)$ has support $\rho$ (so that its projections to the $\rho '$-localizations are zero for any $\rho '\neq\rho$), then its image $f_*(\eta)$ projects to zero in $\k(X,G)_{\sigma '}$ for every $\sigma '\neq\sigma$: indeed for any such $\sigma '$ and $\rho '\twoheadrightarrow \sigma '$ it holds that $\rho '\neq \rho$, and $\eta$ has null image in $\k(Z,  G\times H)_{\rho '}$. We conclude that $\eta$ is zero when localized at $\sigma '$, since $\k(Z,G\times H)_{\sigma '}=\prod_{\substack{\rho ' \in \cC(G\times H)\\
   \rho ' \twoheadrightarrow \sigma '}}\k(Z, G \times H)_{\rho '}$; being $f_*$ a map of $\R G-$modules, $f_*(\eta)$ is also zero when localized at $\sigma '$.

In particular $f_*\eta$ has support at $\sigma$, and thus lies in $\k(X,G)_\sigma$.
\end{enumerate1}

\end{remark}

We will use the following facts, whose proof can be found in \cite{Schadeck}.

\begin{proposition}\label{prop:reduced-pushforward}
Let $i\colon \cX\red \into\cX$ be the reduction of $\cX$. Then
   \[
   i_{*}\colon \kg{\cX\red} \arr \kg\cX
   \]
is an isomorphism.
\end{proposition}

\begin{proposition}\label{prop:descent-geometric-K}
Let $\pi\colon \cX' \arr \cX$ a representable finite faithfully flat morphism of  stacks. Then the sequences
   \[
   0 \arr \kg\cX \xarr{\pi^{*}} \kg{\cX'}
   \xarr{\pr^{*}_{1} - \pr^{*}_{2}}
   \kg{\cX'\times_{\cX} \cX'}   
   \]
and
   \[
   \kg{\cX'\times_{\cX} \cX'}
   \xarr{\pr_{1*} - \pr_{*2}} \kg{\cX'}
   \xarr{\pi_{*}}
   \kg\cX \arr 0   
   \]
are exact.
\end{proposition}

This is easily seen to fail for K-theory, for example, for the representable finite faithfully flat morphism $\spec k \arr \cB_{k}G$, where $G$ is a nontrivial finite group scheme.

\begin{corollary}\label{finite-cover}
Let $\Gamma$ be a finite group, $\pi\colon \cX'\arr \cX$ a Galois cover with group $\Gamma$. Then the pullback
   \[
   \pi^{*} \colon \kg\cX \arr \kg{\cX'}^{\Gamma}
   \]
and the pushforward
   \[
   \pi_{*}\colon (\kg{\cX'})_{\Gamma} \arr \kg\cX
   \]
are isomorphisms.
\end{corollary}

Here of course $\kg{\cX'}^{\Gamma}$ is the group of invariants, and $\pi_{*}(\kg{\cX'})_{\Gamma}$ is the group of covariants.

This applies in particular when $\cX$ is a quotient stack of the form $[X/\Gamma]$, where $\Gamma$ is a finite group acting on a scheme $X$. In this case we have isomorphism $\kg{[X/\Gamma]} \simeq (\k X)^{\Gamma}$ and $(\k X)_{\Gamma} \simeq \kg\cX$.\\

Let $\pi\colon \cX \arr M$ be the moduli space of $\cX$; then we get a homomorphism
   \[
   \pi_{*}\colon \kg\cX \arr \kg M = \k M\,.
   \]

\begin{theorem}\label{thm:pushforward-moduli}
The pushforward $\pi_{*}\colon \kg\cX \arr \k M$ is an isomorphism.
\end{theorem}

\section{The tautological part of the K-theory of the cyclotomic inertia}

Suppose that $K$ is an extension of $k$, and let $(\xi, a)$ be an object of $\ci\cX(K)$. The homomorphism $a\colon \mmu_{\infty, K} \arr \underaut_{K}\xi$ induces an action of $\mmu_{\infty, K}$ on $\xi$, which commutes with itself, since $\mmu_{\infty, K}$ is abelian; thus $a$ can be considered as a homomorphism $\mmu_{\infty, K} \arr \underaut_{K}(\xi, a)$. We can think of this as follows: a morphism $\cB_{K}\mmu_{\infty, K} \arr \cX$ as a canonical lifting to a morphism $\cB_{K}\mmu_{\infty, K} \arr \ci\cX$, which we call its \emph{tautological lifting}.

A morphism $\cB_{K}\mmu_{\infty, K} \arr \ci\cX$ is called \emph{tautological} if it is isomorphic to the tautological lifting of the composite $\cB_{K}\mmu_{\infty, K} \arr \ci\cX \arr \cX$. Equivalently, we can define a tautological morphism as follows. A morphism $\cB_{K}\mmu_{\infty, K} \arr \ci\cX$ corresponds to an object $(\xi, a)$ of $\ci\cX(K)$, together with a homomorphism $b\colon \mmu_{\infty, K} \arr \underaut_{K}(\xi, a)$. This gives two action of $\mmu_{\infty, K}$ on $\xi$: one given by $a\colon \mmu_{\infty, K} \arr \underaut_{K}\xi$, and the other by the composite $\mmu_{\infty, K} \xarr{b} \underaut_{K}(\xi, a) \subseteq \underaut_{K}\xi$. The morphism $\cB_{K}\mmu_{\infty, K} \arr \ci\cX$ is tautological if and only if the two actions coincide.

As we saw, there is an equivalence between morphisms $\cB_{K}\mmu_{\infty, K} \arr \cX$ and representable morphisms $\cB_{K}\mmu_{r, K} \arr \cX$ for some $r$; hence every representable morphism $\cB_{K}\mmu_{r, K} \arr \cX$ has a tautological lifting $\cB_{K}\mmu_{r, K} \arr \ci[r]\cX$. Therefore we can also talk about tautological representable morphisms $\cB_{K}\mmu_{r, K} \arr \ci[r]\cX$.

Let us define two multiplicative systems $\Theta_{\cX} \subseteq \kz(\ci\cX)$ and $\widetilde{\Theta}_{\cX} \subseteq \kz(\ci\cX)$ as follows. An element $\alpha \in \kz(\ci\cX)$ is in $\Theta_{\cX}$ if for every representable tautological morphism $\phi\colon \cB_{K}\mmu_{r, K} \arr \cX$ the image of $\phi^{*}\alpha \in \rR\mmu_{r}$ in $\Q r$ is non-zero. An element $\alpha \in \kz(\ci\cX)$ is in $\widetilde{\Theta}_{\cX}$ if for every non-tautological representable morphism $\phi\colon \cB_{K}\mmu_{r, K} \arr \cX$ the image of $\phi^{*}\alpha \in \rR\mmu_{r}$ in $\Q{r}$ is non-zero. (Here $K$ is an extension of $k$.)

\begin{definition}
The \emph{tautological part} of the K-theory of the cyclotomic inertia stack $\ci\cX$ is the localization $\kt{\ci\cX} \eqdef \Theta_{\cX}^{-1}\k(\ci\cX)$. The \emph{non-tautological part} of the K-theory of the cyclotomic inertia stack $\ci\cX$ is the localization $\knt{\ci\cX} \eqdef {\widetilde{\Theta}_{\cX}}^{-1}\k(\ci\cX)$.
\end{definition}

Alternatively, one can define $\kt{\ci[r]\cX}$ (respectively $\knt{\ci[r]\cX}$) as the localization of $\k(\ci[r]\cX)$ along the multiplicative systems in $\kz{\ci[r]\cX}$ consisting of elements whose image  in $\Q r$ is non-zero for every non-tautological (respectively tautological) representable morphism tautological representable morphisms $\cB_{K}\mmu_{r, K} \arr \ci[r]\cX$. Clearly we have $\kt{\ci\cX} = \bigoplus_{r \geq 1}\kt{\ci[r]\cX}$ and $\knt{\ci\cX} = \bigoplus_{r \geq 1}\knt{\ci[r]\cX}$. Furthermore, if $\cX$ is an algebraic space, then $\kt{\ci \cX} = \k X$ and $\knt{\ci \cX} = 0$.

The following property is also evident.

\begin{proposition}\label{prop:tautological-mu}\hfil
\begin{enumerate1}

\item The projection $\k(\ci[r]\cX) \arr \kt{\ci[r]\cX}$ factors through the projection
   \[
   \k(\ci[r]\cX) \arr \kmu{r}{\ci[r]\cX}\,.
   \]

\item The projection $\k(\ci[r]\cX) \arr \knt{\ci[r]\cX}$ factors through the projection $\k(\ci[r]\cX) \arr \bigoplus_{s \neq r}\kmu{s}{\ci[r]\cX}$.

\end{enumerate1}
\end{proposition}

It is easy to check that the multiplicative systems $\Theta_{\cX} \subseteq \kz(\cI\cX)$ and $\widetilde{\Theta}_{\cX} \subseteq \kz(\cI\cX)$ are invariant under the action of $\autm$ on $\kz$; hence we get an action of $\autm$ on $\kt{\ci\cX}$. By restriction, this gives an action of $\autm$ on each $\kt{\ci[r]\cX}$.

Let us describe the tautological and non-tautological parts of equivariant K-theory of the cyclotomic inertia stack in equivariant terms. Assume that $\cX = [X/G]$, where $G$ is a product of general linear groups. Let $\cU \subseteq \ci[r]\cX$ be a connected component of the cyclotomic inertia stack $\ci\cX$. Recall that $\ci[r]\cX$ is the disjoint union
\[
\ci[r]\cX=\coprod_{\substack {\sigma\in \lc_r(G)\\ |\sigma|=r}}[X^\sigma/C_{G}(\sigma)]
\]
so $\cU$ corresponds to some embedding $\sigma\in\lc_r(G)$. In $\k([X^\sigma/C_{G}(\sigma)])=\k(X,C_{G}(\sigma))$ the multiplicative system $\Theta_\cX$ restricts to the $\Sigma_\sigma$ relative to the embedding $\sigma\into C_{G}(\sigma)$. In particular, from Proposition~\ref{prop:decomposition} we get

\begin{proposition}
If $G$ is a product of general linear groups then 
\[
\kt{\ci[r]\cX}=\bigoplus_{\substack {\sigma\in \overline C(G)}}\k(X^\sigma,C_{G}(\sigma))_\sigma.
\]
\end{proposition}

From this description we obtain the following.

\begin{proposition}
The homomorphism $\k(\ci\cX) \arr \kt{\ci\cX} \times \knt{\ci\cX}$ is an isomorphism.
\end{proposition}

From this splitting we get an embedding $\kt{\ci\cX} \subseteq \kmu r {\ci\cX}$.

%
%
%

\begin{proposition}\label{prop:tame-pushforward}
Let $f\colon \cY \arr \cX$ be a morphism of stacks. 

\begin{enumerate1}

\item If $\cX$ and $\cY$ are regular, the pullback $\ci{f}^{*}\colon \k{(\ci{\cX})} \arr \k{(\ci{\cY})}$ descends to an $\autm$-equivariant homomorphism
   \[
   \ci{f}^{*}\colon \kt{\ci\cX} \arr \kt{\ci\cY}\,.
   \]

\item If $f$ is proper and relatively tame, the pushforward $\ci{f}_{*}\colon \k{(\ci\cY)} \arr \k{(\ci\cX)}$  carries $\kt{\ci\cX} \subseteq \k{(\ci\cX)}$ into $\kt{\ci\cY} \subseteq \k{(\ci{\cY})}$, and induces an $\autm$-equivariant ring homomorphism $\kt{\ci\cX} \arr \kt{\ci\cY}$.

\end{enumerate1}

\end{proposition}

\begin{proof}
The first part follows from the fact that $\ci{f}^*(\Theta_\cX)\subseteq\Theta_\cY$. Let us concentrate now on the second one.

Suppose that $\cY=[Y/H]$ and $\cX=[X/G]$ where $X,Y$ are schemes and $G,H$ are general linear groups.  Recall the following from the proof of Proposition~\ref{prop:functoriality}: there is an $H-$invariant map $Y \arr \cY \xarr{f} \cX$; if $Z = X \times_{\cX}Y$, there is a natural action of $G \times H$ on it, and $[Z/G\times H] = \cY$. The projection $Z \arr X$ is equivariant with respect to $\pi\colon G\times H\arr G$ and induces a map
   \[
   \cY = [Z/ G\times H] \arr [X/G] = \cX
   \]
that is isomorphic to $f$. 

Let $\rho\into G\times H$ be a dual cyclic subgroup, whose projections onto $G,H$ are $\sigma_G,\sigma_H$ respectively. Then $C_{G\times H}(\rho)=C_G(\sigma_G)\times C_H(\sigma_H)$.

Let $\sigma:=\sigma_G$; the map $Z\arr X$ induces a morphism
\[
[Z^\rho/C_G(\sigma_G)\times C_H(\sigma_H)]\arr [X^\sigma/C_G(\sigma)]
\]
and this is exactly the restriction of $\ci{f}$ to the component of $\ci\cY$ relative to $\rho$. The induced map $\ci{f}_*$ is a morphism of $\R C_G(\sigma)$-modules.

We can now argue as in our Observation~\ref{remark:push-forward}: if $\eta\in \k(Z^\rho, C_G(\sigma)\times C_H(\sigma_H))$ has support $\rho$ (so that its projections to the $\rho '$-localizations are zero for any $\rho '\neq\rho$), then its image $\ci{f}_*(\eta)$ has support at $\sigma$, and thus projects to zero in $\knt{[X^\sigma/C_G(\sigma)]}$.

This settles the proof.
\end{proof}

One can also prove that if $f$ is representable, then $\ci{f}_{*}\k{\ci\cY} \arr \k{(\ci\cX)}$ also induces a pushforward $\knt{\ci\cY} \arr \knt{\ci\cX}$, but we will not need this fact.

Let $\cX$ be a stack, and denote by $\rho = \rho_{\cX}\colon \ci{\cX} \arr \cX$ the canonical morphism.

\begin{proposition}\label{lemma:iso}

Let $\cX$ be an algebraic stack with finite cyclotomic inertia. For any $r > 0$, the restriction $\rho_{*}\colon \kt{\ci[r]\cX}^{\autm} \arr \kmu{r}{\cX}$ of the pushforward $\rho_{*}\colon \k(\ci[r]\cX) \arr \k X$ is an isomorphism.

\end{proposition}

\begin{proof}

Let us begin with the following

\begin{lemma}\label{lemma:push-forward}
Let $X$ be a $G-$equivariant scheme and $\sigma\subset G$ a dual cyclic subgroup. Then the push-forward map $(j_\sigma)_*\colon\k(X^\sigma,C_G(\sigma))_\sigma\arr \k(X,C_G(\sigma))_\sigma$ is an isomorphism.
\end{lemma}

\begin{proof}
Let $Y:=X-X^\sigma$; we have an exact sequence
\[
\k[i+1](Y,C_G(\sigma))_\sigma\arr \k[i](X^\sigma,C_G(\sigma))_\sigma\xarr{(j_\sigma)_*} \k[i](X,C_G(\sigma))_\sigma\arr \k[i](Y,C_G(\sigma))_\sigma
\] 
and $Y^\sigma=0$, so $\k[i](Y,C_G(\sigma))_\sigma=0$ for all $i$ by Corollary~\ref{vanishing}{}.
\end{proof}

\step{Step 1} Suppose first that  $\cX=[X/T]$, where $T$ is a split torus and $X$ is a separated scheme. 

By Lemma~\ref{lemma:push-forward}, for any $\sigma\in \cC_r(T)$ the push-forward $(j_\sigma)_*\colon \k(X^\sigma,T)_\sigma\arr \k(X,T)_\sigma$ is an isomorphism. 

But $\rho_*$ factors as
\[
\bigoplus_{\substack {\sigma\in \cC_r(T)}}\k(X^\sigma,T)_\sigma\xarr{\oplus (j_\sigma)_*}\bigoplus_{\substack {\sigma\in \cC_r(T)}}\k(X,T)_\sigma\arr \kmu{r}{X,T}
\]
and the thesis follows.

\step{Step 2} We can finally suppose that $\cX=[X/\GL_{n}]$, where $X$ is a separated scheme; now the naive application of Lemma~\ref{lemma:push-forward} does not give the full result, but we will reduce to the previous case.

Let $T$ be a maximal split subtorus of $\GL_{n}$ and let $\cY:=[X/T]$.

Let us describe the cartesian diagram
\[
   \begin{tikzcd}
   \ci\cX\times_{\cX}\cY \rar["\rho '"]\dar
      & \cY \dar \\
   \ci\cX \rar["\rho"] &\cX
   \end{tikzcd}
   \]
We have
\[
\ci\cX=\coprod_{\substack {\sigma\in \lc(\GL_{n})}}[X^\sigma/C_{\GL_{n}}(\sigma)].
\]
Furthermore any $\sigma\in \lc(\GL_{n})$ has a representative whose image is a subgroup of $T$, whence we may assume that $T\subseteq C_{\GL_{n}}(\sigma)$. Fix any such $\sigma$.

\begin{lemma}
Let $X$ be a $G-$equivariant scheme and $G',G''$ subgroups of $G$. Consider the cartesian diagram
 \[
   \begin{tikzcd}
   \cS\rar["p''"]\dar["p'"]
      & \lbrack X/G'' \rbrack \dar \\
   \lbrack X/G'\rbrack \rar &\lbrack X/G \rbrack
   \end{tikzcd}
   \]
Then $\cS=[(X\times G)/(G'\times G'')]$, with an action given by $(x,g)(\mathbf{g'},\mathbf{g''})=(x\mathbf{g'},\mathbf{g'}^{-1}g\mathbf{g''})$.

The map $p'$ is induced by $(x,g)\arr x$, while $p''$ is induced by $(x,g)\arr xg$.
\end{lemma}
\begin{remark}
 The definition of $\cS$ might look asymmetrical, since  we could make $G'\times G''$ act by $(x,g)(\mathbf{g'},\mathbf{g''})=(x\mathbf{g''},\mathbf{g''}^{-1}g\mathbf{g'})$; however this is not the case. In fact the automorphism of $X\times G$ given by $(x,g)\arr (xg,g^{-1})$ is equivariant with respect to the two actions.
\end{remark}

\begin{proof}
Let us consider the following diagram, where every square is cartesian:
\[
\begin{tikzcd}
X\times G\rar \dar &\cS '\rar\dar & X\dar \\
\cS ''\rar\dar & \cS \rar["p''"]\dar["p'"]&\lbrack X/G'' \rbrack \dar \\
X\rar & \lbrack X/G'\rbrack \rar &\lbrack X/G \rbrack
\end{tikzcd}
\]
The horizontal and vertical maps $X\times G\arr X$ are the multiplication and the projection onto the first factor, respectively.

By base-change, $X\times G\arr \cS'$ is a $G'-$principal bundle. We conclude that $\cS'=[(X\times G)/G']$, where $G'$ acts as $(x,g)\mathbf{g'}=(x\mathbf{g'}, \mathbf{g'}^{-1}g)$.

Analogously, $\cS''=[(X\times G)/G'']$, where $G''$ acts as $(x,g)\mathbf{g''}=(x,g\mathbf{g''})$.

Finally, the map $\cS''\arr\cS$ is a $G'$-principal bundle, so that
\[
\cS=[\cS''/G']=[(X\times G)/(G'\times G'')],
\]
where the action is easily seen to be the one given in the statement.
\end{proof}

Thanks to this lemma, we have a commutative diagram 
\[
   \begin{tikzcd}
    \lbrack(X^\sigma\times \GL_n)/(C_{\GL_n}(\sigma)\times T)\rbrack\rar \dar["\pi '"] &
\lbrack(X\times \GL_n)/(C_{\GL_n}(\sigma)\times T)\rbrack\rar \dar["\pi"] &
\lbrack X/T\rbrack \dar["\pi"] \\
   \lbrack X^\sigma /C_{\GL_n}(\sigma)\rbrack \rar & \lbrack X/C_{\GL_n}(\sigma)\rbrack \rar &\lbrack X/\GL_n\rbrack
   \end{tikzcd}
   \]
where the squares are cartesian.

Let us consider the conjugacy classes of dual cyclic subgroups $\tilde\sigma\subset (C_{\GL_n}(\sigma)\times T)$ who lie over $\sigma$ and have fixed points in $X^\sigma\times\GL_n$. If $\tilde\sigma$ is such an element, let $G'\into \GL_n$ be the projection onto $\GL_n$ of the fixed locus of $\tilde\sigma$. Calling $\sigma '$ the projection of $\tilde\sigma$ on $T$, we must have that the conjugation map $G'\times \GL_n\arr\GL_n$, $(g,s)\arr g^{-1}sg$, sends $\sigma$ to $\sigma '$.

We conclude that the sought $\tilde\sigma$ are exactly those given by the embeddings $\sigma\arr (\sigma,\sigma ')\subset (C_{\GL_n}(\sigma)\times T)$, where $\sigma '$ is a conjugate of $\sigma$ lying in $T$. Besides, we can assume that $\sigma '=g^{-1}\sigma g$, where $g\in S_n$ is a permutation matrix, thanks to the following Observation~\ref{lemma:centralizer}:

\begin{observation}\label{lemma:centralizer}
$\cC(G)$ is in natural bijective correspondence with the set of orbits for the action of $\rS_{n}$ on $\cC(T)$.
\end{observation}

The Weyl group $W=S_n=N_{\GL_n}(T)/T$ acts on the set of the $\tilde\sigma$'s, and the stabilizer of any element is the Weyl group of $C_{\GL_n}(\sigma)$, which we will denote as $\Delta_\sigma$.

For any such $\tilde\sigma=(\sigma,g^{-1}\sigma g)$ there is a commuting diagram
\[
   \begin{tikzcd}
\lbrack X^\sigma/T\rbrack\rar["i_{\tilde\sigma}"] \arrow[rd] & \lbrack(X^\sigma\times \GL_n)/(C_{\GL_n}(\sigma)\times T)\rbrack \dar["\pi"] \\
&\lbrack X^\sigma /C_{\GL_n}(\sigma)\rbrack .
 \end{tikzcd}
   \]
There is a group morphism $T\arr C_{\GL_n}(\sigma)$, $t\arr (t,g^{-1}tg)$, and $i_{\tilde\sigma}$ is induced by $x\arr (x,g)$, which is equivariant with respect to the previous map.  

The subscheme of $\tilde\sigma$-fixed points in $X^\sigma\times \GL_n$ is $X^\sigma\times C_{\GL_n}(\sigma)g$ and $\tilde\sigma$ is central in $C_{\GL_n}(\sigma)\times T$. So the map
\[
[X^\sigma/T]=[X^\sigma\times C_{\GL_n}(\sigma)/C_{\GL_n}(\sigma)\times T]\arr [X^\sigma\times\GL_n/C_{\GL_n}(\sigma)\times T]
\]
is just the morphism induced from the inertia map
\[
\rho\colon \ci {[X^\sigma\times\GL_n/C_{\GL_n}(\sigma)\times T]}\arr [X^\sigma\times\GL_n/C_{\GL_n}(\sigma)\times T].
\]

Let us call, again, $\sigma '$ the dual cyclic subgroup $g^{-1}\sigma g\subset T$. By Lemma~\ref{lemma:push-forward}, the maps
\[
(j_{\tilde\sigma})_*=(\rho '\circ i_{\tilde\sigma})_*\colon \k(X^\sigma,T)_{\sigma '}\arr \k(X,T)_{\sigma '}
\]
and
\[
(i_{\tilde\sigma})_*\colon \k(X^\sigma,T)_{\sigma '}\arr \k(X^\sigma\times \GL_n,C_{\GL_n}(\sigma)\times T)_{\tilde\sigma}
\]
are isomorphisms.

In particular $\rho '_*\colon \k(X^\sigma\times \GL_n,C_{\GL_n}(\sigma)\times T)_{\tilde\sigma}\arr \k(X,T)_{\sigma '}$ is an isomorphism.\\

Since the map $\pi^*$ induces an isomorphism $\k(X^\sigma, C_{\GL_n}(\sigma))\simeq\k(X^\sigma\times\GL_n, C_{\GL_n}(\sigma)\times T)^W$, we conclude by Observation~\ref{remark:push-forward} that for any $\tilde\sigma$ above $\sigma $ it gives an isomorphism
\[
\k(X^\sigma, C_{\GL_n}(\sigma))_\sigma\simeq\k(X^\sigma\times\GL_n, C_{\GL_n}(\sigma)\times T)_{\tilde\sigma}^{\Delta_\sigma}
\]
Let $\Gamma_\sigma=(N_{\GL_n}(T)\cap N_{\GL_n}(\sigma))/T$ be the stabilizer of $\tilde\sigma\subset T$ for the action of $W$ on $\mathcal C(T)$. Again by Observation~\ref{remark:push-forward} we infer that $\pi^*$ induces an isomorphism (see Appendix B)
\[
\k(X, \GL_n)_\sigma\simeq \k(X,T)_{\sigma '}^{\Gamma_\sigma}
\]
Now, by Lemma~\ref{lemma:centralizer} there is an exact sequence
\[
0\arr \Delta_\sigma\arr \Gamma_\sigma\arr w_{\GL_n}(\sigma)\arr 0
\]
so we have a commuting diagram (as flat pull-backs and proper push-forwards commute in a cartesian square)
\[
\begin{tikzcd}
\k(X^\sigma, C_{\GL_n}(\sigma))_\sigma^{w_{\GL_n}(\sigma)} \rar["\rho_*"]\dar["\pi^*"] & \k(X,\GL_n)_\sigma \dar["\pi^*"] \\
\k(X^\sigma\times \GL_n,C_{\GL_n}(\sigma)\times T)_{\tilde\sigma}^{\Gamma_\sigma}\rar["\rho '_*"] & \k(X,T)_{\sigma '}^{\Gamma_\sigma} 
\end{tikzcd}
\]
The two vertical arrows are isomorphisms and so is $\rho '_*$. We see then that $\rho_*$ is also an isomorphism, and this concludes the proof.

\end{proof}

\begin{corollary}
The pushforward $\rho_{*}$ induces a group isomorphism
   \[
   \rho_{*}\colon \kt{\ci\cX}^{\autm} \simeq \k\cX\,.
   \]

Furthermore, $\rho_{*}$ is covariant for proper morphisms.
\end{corollary}

\section{The twist map}

If $\cX$ is a noetherian algebraic stack, we denote by $\coh \cX$ the category of coherent sheaves on $\cX$. If $G$ is a finite diagonalizable group scheme over $R$ with character group $\widehat{G} = \hom(G, \gm)$, then the K-theory group of the product $\cX \times \cB_{R} G$ is the tensor product $\k{(\cX)} \otimes_{\k(A)} \R G$. The point here is that the category $\coh{(\cX \times \cB_{k} G)}$ is the category of coherent sheaves on $\cX$ with an action of $G$, which in turn can be identified with the direct sum of categories $\oplus_{\chi \in \widehat{G}} \coh\cX$, by the usual splitting of quasi-coherent sheaves with an action of $G$ into eigensheaves. If $F$ is a coherent sheaf on $\cX$ and $\chi\in \widehat{G}$, we denote by $F\otimes \chi$ the sheaf $F$ with the action of $G$ defined by $\chi$, considered as a sheaf on $\cX \otimes \cB_{R}G$.

Fix a positive integer $r$. There is a morphism
   \[
   \alpha_{\cX}\colon \ci[r]\cX \times \cB_{R}\mmu_{r} \arr \ci[r]\cX
   \]
that is defined, at the level of objects, as follows. Let $(\xi, a, P)$ be an object of $\ci[r]\cX \times \cB_{R}\mmu_{r}$ over an affine scheme $S$; here $\xi$ is an object of $\cX(S)$, $a\colon \mmu_{r, S} \arr \underaut_{S}\xi$ is a monomorphism of group schemes over $S$, and $P \arr S$ is a $\mmu_{r}$-torsor. By descent along $P \arr S$ we obtain another object $\xi_{P}$ of $\cX(S)$, which we can think about as the quotient $(\xi \times P)/\mmu_{r}$; the automorphism group scheme $\underaut_{S}\xi_{P} \arr S$ is the twisted form $(\underaut_{S}\xi)_{P}$ obtained by descent along $P \arr S$ using the action of $\mmu_{r, S}$ on $\underaut_{S}\xi$ by conjugation. Since this action by conjugation leaves $a$ invariant, we obtain a homomorphism $a_{P}\colon \mmu_{r, S} \arr \underaut_{S}\xi_{P}$. This gives an object $\alpha_{\cX}(\xi, a, P) = (\xi_{P}, a_{P}, P)$. The effect of $\alpha_{\cX}$ on arrows is clear.\\
We will refer to this morphism as the \emph{twist map}.\\

Equivariantly, if $\cX=[X/G]$, on any component $[X^\sigma/C_G(\sigma)]\times\cB_R\sigma$ the map $\alpha_\cX$ is given as follows: there is a projection
\[
m\colon C_G(\sigma)\times\sigma\arr C_G(\sigma),
\]
corresponding to the multiplication map (this is a morphism of groups), which induces a map $[X^\sigma/C_G(\sigma)\times\sigma]\arr[X^\sigma/C_G(\sigma)]$. This coincides with $\alpha_\cX$.\\

Now consider the disjoint union 

\[
\widetilde{\ci\cX} = \coprod_{\substack{r}}\ci[r]\cX\times\cB_R\mmu_r.
\]
For any $f\colon\cY\arr\cX$ the map $\ci{f}$ induces a morphism $\ci\cY\times\cB_R\mmu_{\infty,R}\arr \ci\cX\times\cB_R\mmu_{\infty,R}$ and by projection a map
\[
\widetilde{\ci{f}}\colon \widetilde{\ci\cY}\arr \widetilde{\ci\cX}.
\]
Let us define a multiplicative system $\Omega_\cX\subseteq \k[0](\widetilde{\ci\cX})$ as follows: for any component $\ci[r]\cX\times\cB_R\mmu_r$ consider the representable morphisms $\phi$ from $ \cB_R\mmu_r$ into it, that are trivial on $\ci[r]\cX$ and the identity on $\cB_R\mmu_r$; $\alpha\in \k[0](\ci[r]\cX\times\cB_R\mmu_r)$ lies in $\Omega_\cX$ if and only if $\phi^*\alpha$ has a nontrivial projection on $\Rt\mmu_r$ for any such $\phi$. $\Omega_\cX$ is the multiplicative systems generated by those elements for any $r$.

The localization $\Omega_\cX^{-1}\k(\widetilde{\ci\cX})$ can be described in equivariant terms. If $\cX=[X/G]$ ($G$ a product of general linear groups) and $\sigma\subset G$ is a dual cyclic subgroup then $\Omega_\cX^{-1}\k([X^\sigma/C_G(\sigma)]\times\cB_R\mmu_r)$ coincides with the localization - as a $\R(C_G(\sigma)\times \mmu_r)-$module at the dual cyclic subgroup $\tilde\sigma=\text{Id}\times\mmu_r\into C_G(\sigma)\times\mmu_r$.

\begin{definition}
Let the \emph{tautological part} $\kt{\widetilde{\ci\cX}}$ of $\k(\widetilde{\ci\cX})$ be the localization $\Omega_\cX^{-1}\k(\widetilde{\ci\cX})$.
\end{definition}

The projection $\k(\widetilde{\ci\cX})\arr \kt{\widetilde{\ci\cX}}$ has a natural splitting, which we shall describe below.

For any $d|r$ we have that $\kmu d {\mmu_r}\simeq \k(A)\otimes\QQ(\zeta_d)$, hence there exists an immersion 
\[
i_d:\QQ(\zeta_d)\into \R(\mmu_r).
\]
We can describe it explicitly: let us define the polynomial
\[
\psi_d:=\frac{t^r-1}{\phi_d}=\prod_{\substack{d_1|r \\ d_1\neq d}}\phi_{d_1}
\]
and note that $\psi_d(\zeta_d)\neq 0$; then for any $x\in\QQ(\zeta_d)$ we can choose a polynomial $p_x$ such that $p_x(\zeta_d)=x\cdot \psi_d(\zeta_d)^{-1}$ and it holds that
\[
i_d(x)=p_x\cdot\psi_d(t) \in \QQ[t]/(t^r-1)= \R(\mmu_r)
\]
This definition is manifestly independent of the choice of $p_x$, since it is well-defined up to multiples of $\phi_d(t)$ and $\psi_d\cdot \phi_d(t)=0$ in $\R(\mmu_r)$.

In particular $i_1\colon\QQ \arr \QQ[t]/(t^r-1)$ sends $q\in\QQ$ to $q\frac{1+t+\dots +t^{r-1}}{r}$.\\

Since we also have an immersion $\kg{\ci[r]\cX}\into \k(\ci[r]\cX)$ we can combine it with $i_r$ to get a map
\[
\kg{\ci[r]\cX}\otimes\QQ(\zeta_r)\into \k(\ci[r]\cX)\otimes_{K_0(A)} \R(\mmu_r)=\k(\ci[r]\cX \times \cB_{R}\mmu_{r}).
\]
This is the required splitting, since $S_{\tilde\sigma}^{-1}\k(\widetilde{\ci[r]\cX})=\kg{\ci[r]\cX}\otimes\QQ(\zeta_r)$.\\

\begin{proposition}
Let $f\colon\cY\arr\cX$ be a proper map of quotient stacks. Then $\widetilde{\ci{f}_*}$ carries $\kt{\widetilde{\ci\cY}}$ into $\kt{\widetilde{\ci\cX}}$ and induces a ring homomorphism on the tautological parts.
\end{proposition}

\begin{proof}
As we already did, we can suppose that $\cY=[Y/H]$ and $\cX=[X/G]$ where $X,Y$ are schemes and $G,H$ are general linear groups.  If $Z = X \times_{\cX}Y$, there is an action of $G \times H$ on it, and $[Z/G\times H] = \cY$. The projection $Z \arr X$ is equivariant with respect to $\pi\colon G\times H\arr G$ and induces a map
   \[
   \cY = [Z/ G\times H] \arr [X/G] = \cX
   \]
that is isomorphic to $f$. 

Let $\rho\into G\times H$ be a dual cyclic subgroup, whose projections onto $G,H$ are $\sigma_G,\sigma_H$ respectively. Then $C_{G\times H}(\rho)=C_G(\sigma_G)\times C_H(\sigma_H)$.

Let $\sigma:=\sigma_G$; the maps $Z\arr X$ and $\cB_R\rho\arr\cB_R\sigma$ induce a morphism
\[
[Z^\rho/C_G(\sigma_G)\times C_H(\sigma_H)]\times \cB_R\rho\arr [X^\sigma/C_G(\sigma)]\times\cB_R\sigma
\]
and this is exactly the restriction of $\widetilde{\ci{f}}$ to the component of $\widetilde{\ci\cY}$ relative to $\rho$. The induced map $\widetilde{\ci{f}_*}$ is a morphism of $\R( C_G(\sigma)\times\sigma)$-modules.

Suppose that $\nu\in\k([Z^\rho/C_G(\sigma_G)\times C_H(\sigma_H)]\times \cB_k\rho)$ has support at $\tilde\rho=\text{id}\times\rho\subset (C_G(\sigma_G)\times C_H(\sigma_H))\times\rho$. Then $\widetilde{\ci{f}_*}(\nu)$ has support at $\tilde\sigma$($=\text{id}\times\sigma$), by the argument of Observation~\ref{remark:push-forward}.

\end{proof}

Composing the covariant inclusion $\kt{\widetilde{\ci\cX}}\into \k(\widetilde{\ci\cX})$ with
$(\alpha_{\cX})_*\colon \k(\widetilde{\ci\cX}) \arr \k(\ci\cX)$
we get a map 

\[
\beta_\cX\colon \bigoplus_{\substack{r}}\kg{\ci[r]\cX}\otimes\QQ(\zeta_r)\arr \k(\ci\cX).
\]
This is clearly covariant with respect to proper push-forwards.

\section{The main result}

We are now ready to state the main results of this paper, that is a Riemann-Roch formula for the $K-$theory of quotient stacks with finite cyclotomic inertia.

Let $\cX$ be such a quotient stack. Then we have a map given by the composition
\[
\k(\cX)\xrightarrow{\rho^*} \kt{\cI\cX}^{\autm}\xrightarrow{\alpha_\cX^*}\Bigl(\bigoplus_{\substack{r}}\kg{\ci[r]\cX}\otimes\QQ(\zeta_r) \Bigr)^{\autm}
\]
(the last arrow is the composition of $\alpha_\cX^*$ and the projections $\R\mmu_r\arr\QQ(\zeta_r)$). The main result of \cite{vezzosi-vistoli02} says that this is a contravariant isomorphism of algebras. Our goal is to replace it with a covariant isomorphism of $\QQ-$modules, which we will call the
 \emph{Lefschetz-Riemann-Roch} isomorphism.
We will prove the following

\begin{proposition}\label{twist}
$\beta_\cX$ factors through $\kt{\ci\cX}\into \k(\ci\cX)$ and induces an isomorphism
\[
\beta_\cX\colon \bigoplus_{\substack{r}}\kg{\ci[r]\cX}\otimes\QQ(\zeta_r) \arr \kt{\ci\cX}.
\]
\end{proposition}

Recall from Proposition~\ref{lemma:iso} that there is a group isomorphism
\[
\rho_{*}\colon \kt{\ci\cX}^{\autm} \simeq \k(\cX)
\]
so we have a covariant map
\[
\rho_*^{-1}\colon \k(\cX)\arr \kt{\ci\cX}^{\autm}. 
\]

Let $p_*\colon \Bigl(\bigoplus_{\substack{r}}\kg{\ci[r]\cX}\otimes\QQ(\zeta_r) \Bigr)^{\autm}\arr \k(\cX) $ be the isomorphism given by the composition $\rho_*\circ \beta_{\cX}$. Then the map $\cL :=p_*^{-1}$ gives

\begin{theorem}\label{main}
Let $\cX$ be a quotient stack with finite cyclotomic inertia over a base scheme $A=Spec(R)$. Then the above map gives a
 \emph{Lefschetz-Riemann-Roch} isomorphism, which is covariant with respect to proper push-forwards of relatively tame morphisms:
\[
\cL\colon \k(\cX)\longrightarrow \Bigl(\bigoplus_{\substack{r}}\kg{\ci[r]\cX}\otimes\QQ(\zeta_r) \Bigr)^{\autm}
\]
\end{theorem}

The formula $\cL ^{-1}= \rho_*\circ \beta_{\cX}$ is, of course, valid for any quotient stack. However if $\cX$ is regular we can express $\cL$ in a different way, much more feasible for calculations.

Let us consider the map $p^*$ given by the composite
\[
\k(\cX)\xrightarrow{\lambda_{-1}(\cN)^{-1}\cdot\rho^*} \kt{\cI\cX}^{\autm}\xrightarrow{\alpha_\cX^*}\Bigl(\bigoplus_{\substack{r}}\kg{\ci[r]\cX}\otimes\QQ(\zeta_r) \Bigr)^{\autm}
\]
where $\cN$ is the conormal sheaf for $\rho\colon\cI\cX\arr\cX$ and the last map is - obviously - composed with the projections $\R(\mmu_r)\arr\QQ(\zeta_r)$.\\

Our comparison theorem is as follows:

\begin{theorem}\label{comp}
Let $\cX$ be a regular stack with finite regular cyclotomic inertia over $A$. Suppose that $\cX$ is either
\begin{enumerate}
\item of the form $[X/G]$, where $X$ is a scheme and $G$ is a finite flat group scheme over $A$, such that all $\sigma$'s come from subgroups of $G$,
\item Deligne-Mumford
\item tame.

\end{enumerate}  

Then the composition $$p^*\circ p_*\colon \Bigl(\bigoplus_{\substack{r}}\kg{\ci[r]\cX}\otimes\QQ(\zeta_r) \Bigr)^{\autm}\arr \Bigl(\bigoplus_{\substack{r}}\kg{\ci[r]\cX}\otimes\QQ(\zeta_r) \Bigr)^{\autm}$$
is equal to the endomorphism that on each component $\kg{\ci[r]\cX}\otimes\QQ(\zeta_r)$ is a multiplication by the rational number $\frac{\phi(r)}{r}$.
\end{theorem}

In particular the Lefschetz-Riemann-Roch map, in the regular case, is equal to $\bigoplus_{\substack{d}}p^*\cdot \frac{d}{\phi(d)}$.\\

\begin{remark}
If $\cX$ is regular, then the cyclotomic inertia is also regular. Indeed, taking $\cX=[X/\GL_n]$ we have that $X$ is regular and, by \cite{Iversen}, Lemma $2.3$, $X^\sigma$ is also regular for any dual cyclic subgroup $\sigma\subset \GL_n$.
\end{remark}

We are now ready to give the proof of our main result, which amounts to the proof of Proposition~\ref{twist}:  

\begin{proof}
We can suppose $\cX=[X/G]$, where $X$ is a scheme and $G$ is an algebraic group. 

Recall that $\kg{X^\sigma, C_G(\sigma)}\otimes\QQ(\zeta_r)$ is the localization of $\k(X^\sigma, C_G(\sigma)\times\sigma)$ at $\tilde\sigma=\text{id}\times\sigma$. The image of $\tilde\sigma$ for the multiplication map $m\colon C_G(\sigma)\times\sigma\arr C_G(\sigma) $ is $\sigma$ so that, by Observation~\ref{remark:push-forward},  $\beta_\cX$  induces for any $\sigma\in \lc_r(G)$ a map
\[
\beta_{\cX,\sigma}\colon \kg{X^\sigma, C_G(\sigma)}\otimes\QQ(\zeta_r) \arr \k(X^\sigma, C_G(\sigma))_\sigma.
\]
We are left to show that these are isomorphisms.

\step{Step 1} We treat the case where $\cX$ is the classifying space of a finite diagonalizable group scheme. For any $i$ we have $\k[i](\cX)=\k[i](A)\otimes_{\k[0](A)}\k[0](\cX)$, so we just need to prove the theorem for $i=0$.\\

Suppose first that $\cX=\cB_{R}\mmu_n$, for any $n$. In this case we can just verify the claim with a very easy calculation.

For any $d|n$ the map $\beta_{\cX,\langle d\rangle}$ does the following: firstly we have an immersion
\[
\kg[0]{\mmu_n}\otimes\QQ(\zeta_d)=\QQ\otimes\QQ(\zeta_d)\into \R(\mmu_n)\otimes\R(\mmu_d)=\QQ[t]/(t^n-1)\otimes\QQ[s]/(s^d-1)
\]
which sends $q\otimes x$ to $q\cdot\bigl( \frac{1+t+\dots+t^{n-1}}{n}\bigr)\otimes (p_x\cdot\psi_d(s))$.

Finally we compose this map with $(\alpha_\cX)_*\colon \R(\mmu_n)\otimes\R(\mmu_d)\arr \R(\mmu_n)$.

In this case we have a very explicit characterization of $(\alpha_\cX)_*$: indeed it is the push-forward induced by the multiplication morphism
\[
m_d\colon \mmu_d\times \mmu_n\arr \mmu_n
\]
and it just takes the invariants relative to the subgroup $\ker(m_d)=\{(x,x^{-1})\colon x\in\mmu_d\}$. In particular it is quite easy to convince ourself that $(\alpha_\cX)_*$ sends $t^i\otimes s^j$ to $t^i$ if $i\equiv j\pmod{d}$ and to $0$ otherwise.

We see then that the restriction
\[
(\alpha_\cX)_*\colon \R(\mmu_d)\simeq \Bigl( \frac{1+t+\dots+t^{n-1}}{n}\Bigr)\otimes \R(\mmu_d)\arr \R(\mmu_n)
\]
sends the image of a polynomial $p$ to 
\[
p\cdot\frac{1}{n}(1+t^{d}+\dots+t^{n-d})=p\cdot\frac{t^n-1}{n(t^d-1)}
\]
which is nothing but $1/n$ times the push-forward $\R(\mmu_d)\arr \R(\mmu_n)$ induced by the immersion $\mmu_d\into \mmu_n$.

Now consider the map $\cB\mmu_d\arr\cB \mmu_n$ induced by the immersion $\mmu_d\into \mmu_n$; it is a finite representable map, and induces isomorphisms $\QQ\simeq (\R\mmu_d)_\textbf{g}\arr (\R \mmu_n)_\textbf{g}\simeq\QQ$ and $(\R\mmu_d)_{\mmu_d}\arr (\R\mmu_n)_{\mmu_d}$.
We thus get a commutative diagram
\[
\begin{tikzcd}
\QQ(\zeta_d)\otimes (\R\mmu_d)_\textbf{g}\rar\dar & \R\mmu_d\otimes \R \mmu_d\rar["\beta_{\cB\mmu_d}"]\dar & (\R \mmu_d)_{\mmu_d}\dar  \\
\QQ(\zeta_d)\otimes (\R \mmu_n)_\textbf{g}\rar & \R\mmu_d\otimes \R \mmu_n\rar["\beta_{\cB \mmu_n}"] & (\R \mmu_n)_{\mmu_d} 
\end{tikzcd}
\]

By the above computation, the top row is equal to $\frac{1}{d}\colon\QQ(\zeta_d)\arr\QQ(\zeta_d)$, and is an isomorphism.   We conclude that $\beta_{\cX,\langle d\rangle}\colon \QQ(\zeta_d)\arr \R(\mmu_n)_{\mmu_d}$ is also an isomorphism.\\

Suppose now that $\cX=\cB_R \mathbf{\Delta}$ is the classifying stack of a finite diagonalizable group scheme of order $n$.

Let us fix an immersion $i\colon\mmu_d\arr\mathbf \Delta$. The injection $\QQ=\R(\mathbf\Delta)_\text{\textbf g}\into \R(\mathbf \Delta)$ sends $1$ to $1/n$ times the regular representation of $\mathbf \Delta$: this follows immediately from the commutativity of the diagram

\[
   \begin{tikzcd}
   \QQ=\kg[0]{\cB_R\mathbf{1}} \rar["j_*"]\dar
      & \QQ=\kg[0]{\cB_R\mathbf \Delta} \dar \\
   \QQ=\k[0](\cB_R\textbf 1) \rar["j_*"] &\k[0](\cB_R\mathbf \Delta)
   \end{tikzcd}
   \]
where $j\colon\cB_R\textbf 1\arr \cB_R\mathbf\Delta$ is the map induced by the inclusion of the identity subgroup. Indeed $j^*j_*\colon \kg[0]{\cB_R\mathbf{1}}\arr \kg[0]{\cB_R\mathbf{1}}$ is the multiplication by $n$ map and $j^*$ is an isomorphism of fields, so the above $j_*$ is a multiplication by $n$. 

Then we conclude, as before, that the restriction
\[
(\alpha_\cX)_*\colon \R(\mmu_d)\simeq \Bigl( \frac{\sum_{\delta\in\chi(\mathbf \Delta)}\delta}{n}\Bigr)\otimes \R(\mmu_d)\arr \R(\mathbf \Delta)
\]
is $1/n$ times the push-forward map $i_*$ induced by $i\colon\mmu_d\arr\mathbf\Delta$: indeed $(\alpha_\cX)_*(\delta\otimes s^i)=\delta$ if $\delta_{|\mmu_d}=s^i$ and is $0$ otherwise.

We then see, exactly as before, that in the general case $\beta_{\cX,\mmu_d}\colon \QQ(\zeta_d)\arr \R(\mathbf\Delta)$ is still a multiple of the immersion induced by the canonical decomposition of $R(\mathbf\Delta)$ and thus an isomorphism into the $\mmu_d$-part.\\

\step{Step 2} Now we can handle the case where $\cX=[X/T]$ is the quotient of an algebraic variety by a split torus. We proceed by noetherian induction on $X^\sigma$, the base case where $X^\sigma$ is empty being tautological.

By Thomason's generic slice theorem for torus actions there exists a $T$-invariant nonempty open subscheme $U\subset X^\sigma$, a subgroup $T'$ of $T$ and a $T$-equivariant isomorphism 
\[
U^\sigma=U\simeq T/T'\times (U/T).
\]
Since $T$ acts on $X$ with finite stabilizers $T'$ is finite and diagonalizable. Let $Z=Z^\sigma:=X^\sigma \backslash U$; there is a localization exact sequence

 \[
   \begin{tikzcd}[column sep= 17pt]
    \ar[r]&
     \kg[i]{Z,T}\otimes\Rt(\sigma)  \ar[r]\ar[d, "\beta_{Z,\sigma}"]
      & \kg[i]{X,T}\otimes\Rt(\sigma) \ar[r]\ar[d, "\beta_{X,\sigma}"]
      & \kg[i]{U,T}\otimes\Rt(\sigma) \ar[r]\ar[d, "\beta_{U,\sigma}"]
      &\textblank\\
   \ar[r]&\k[i](Z,T)_\sigma \ar[r]& \k[i](X, T)_\sigma \ar[r]&\k[i](U, T)_\sigma \ar[r]&\textblank
   \end{tikzcd}
   \]

By the five lemma and the inductive hypothesis we just need to prove the theorem for $U$. But in this case $\k(U,T)\simeq \kg{U/T}\otimes \R(T')$ and the map $\beta_{\cU,\sigma}$ is induced by 
\[
\beta_{\cB_R({T'}),\sigma}\colon \R(T')_{\bfg}\otimes \Rt(\sigma)\arr \R(T')_\sigma
\]
We have reduced ourselves to the case of the first step, which has already been handled.

\step{Step 3} We can now conclude the proof dealing with the case $\cX=[X/\GL_n]$: it follows almost immediately from the previous one.

Let $T\subset\GL_n$ be a maximal split torus and $\cY:=X/T$; as always we can assume that $\sigma$ factors through $T$. Let $\pi\colon X^\sigma/C_{\GL_n}(\sigma)\arr X^\sigma/T$ be the projection; this map is representable, so $\pi^*$ is a map of $\R(\GL_n)$-modules and preserves the fundamental decomposition of $\k(X^\sigma, C_{\GL_n}(\sigma))$.

Let us be more precise. Let $\Delta_\sigma$ be the Weyl group of $C_{\GL_n}(\sigma)$; any $\sigma '\in C(C_{\GL_n}(\sigma))$ can be identified with an orbit of the action of $\Delta_\sigma$ on $C(T)$: we can pick any representative for the elements of this orbit, and with an abuse of notation we will call it $\sigma '$ again. Let $\Gamma_{\sigma '}\subseteq\Delta_\sigma$ be the stabilizer of $\sigma '$; the map $\pi^*$ sends $\k(X^\sigma, T)_{\sigma '}^{\Gamma_{\sigma '}}$ isomorphically to $\k(X^\sigma,C_{\GL_n}(\sigma))_{\sigma '}$.\\

Obviously $\Gamma_{\mathbf 1}=\Gamma_\sigma=\Delta_\sigma$ and there is a commutative diagram
\[
   \begin{tikzcd}
   \kg{X^\sigma, C_G(\sigma)}\otimes\QQ(\zeta_r) \rar["\beta_{\cX,\sigma}"]\dar["\pi^*"]
      & \k(X^\sigma, C_G(\sigma))_\sigma \dar["\pi^*"] \\
   \kg{X^\sigma, T}^{\Delta_{\sigma }}\otimes\QQ(\zeta_r) \rar["\beta_{\cY,\sigma}"] & \k(X^\sigma, T)_{\sigma }^{\Delta_{\sigma }}
   \end{tikzcd}
   \]

To conclude the proof we just need to see that the bottom map is an isomorphism: this follows from the previous case.
\end{proof}

Now we prove Theorem~\ref{comp}. It is an immediate consequence of the two Lemmas we give below.\\

\begin{lemma}
The composition
\[
\kt{\ci[r]\cX}^{\autm}\xrightarrow{\rho_*}\k(\cX)\xrightarrow{\lambda_{-1}(\cN)^{-1}\cdot\rho^*} \kt{\ci[r]\cX}^{\autm}
\]
is a multiplication by $\phi(r)$.
\end{lemma}

\begin{remark}
This result is nontrivial, even if we assume the self-intersection formula. Consider for example the case $\cX=\cB S_3$ and $r=2$; we have $\ci[2]\cB S_3=\cB\mmu_2$ and $\rho^*\rho_*\colon \R\mmu_2\arr\R\mmu_2$ sends $f\in\QQ[t]/(t^2-1)$ to $f+f(1)(1+t)$. Thus the result only holds, in general, after passing to the tautological part.
\end{remark}

\begin{proof}
Consider first the case where $\cX=[X/T]$ is the quotient of a scheme by the action of a torus. Let us fix a dual cyclic subgroup $\sigma\simeq\mmu_r\subset T$ and consider the $\sigma-$local part of the k-theory $\k(\cX)$. We have a map
\[
\Bigl(\underset{\phi\in\text{Aut}(\mmu_r)}{\bigoplus}\k(X^{\phi(\sigma)},T)_\sigma\Bigr)^{\text{Aut}(\mmu_r)}\xrightarrow{\rho_*}\k(X,T)_\sigma
\]
and, obviously, $\text{Aut}(\mmu_r)$ acts freely on the left factors. In particular an invariant element is of the form $(x,x,\dots)$ (where we have identified the factors via the obvious isomorphisms) and its image in $\k([X/T])$ is $\phi(r)\cdot i_*(x)$, where $i\colon [X^\sigma/T]\into [X/T]$ is a closed embedding. By the self-intersection formula we have 
\[
\lambda_{-1}(\cN)^{-1}\cdot i^*i_*(x)=x
\]
so the composition
\[
\Bigl(\underset{\phi\in\text{Aut}(\mmu_r)}{\bigoplus}\k(X^{\phi(\sigma)},T)_\sigma\Bigr)^{\text{Aut}(\mmu_r)}\xrightarrow{\rho_*}\k(X,T)_\sigma\xrightarrow{\lambda_{-1}(\cN)^{-1}\cdot\rho^*}
\Bigl(\underset{\phi\in\text{Aut}(\mmu_r)}{\bigoplus}\k(X^{\phi(\sigma)},T)_\sigma\Bigr)^{\text{Aut}(\mmu_r)}
\]
sends $(x,x,\dots)$ to $(\phi(r)x,\phi(r)x,\dots)$, as we wanted to see.\\

Now let us deal with the general case, and suppose $\cX=[X/\GL_n]$. Let us fix a split maximal torus $T\subset\GL_n$ and a dual cyclic subgroup $\sigma\simeq\mmu_r\subset T$. We saw that there are commutative diagrams
\[
   \begin{tikzcd}
    \lbrack X^\sigma/ T \rbrack\arrow[rd,swap,"\pi_\sigma"]\rar &
\lbrack(X^\sigma\times \GL_n)/(C_{\GL_n}(\sigma)\times T)\rbrack\rar["\rho '"] \dar["\pi '"] &
\lbrack X/T\rbrack \dar["\pi"] \\
   & \lbrack X^\sigma/C_{\GL_n}(\sigma)\rbrack \rar["\rho"] &\lbrack X/\GL_n\rbrack
   \end{tikzcd}
\]
where the right square is cartesian. The action of $C_{\GL_n}(\sigma)\times T$ on $X^\sigma\times \GL_n$ is $(\mathbf c,\mathbf t)(x,g)=(x\mathbf c,\mathbf c^{-1}g\mathbf t)$, $\pi '$ is induced by the projection on the first factor and $\rho '$ is induced by $(x,g)\arr x\cdot g$.

For any $g\in\GL_n$ permutation matrix we can form a dual cyclic subgroup $\tilde\sigma =(\sigma,\sigma ')=(\sigma,g^{-1}\sigma g)\subset C_{\GL_n}(\sigma)\times T$ and these are the only ones who lie over $\sigma$. Moreover, for any such $g$, there is a homomorphism $T\arr C_{\GL_n}(\sigma)\times T $, $t\arr (t,g^{-1}tg)$ and a closed embedding $i_{\tilde\sigma}\colon [X^\sigma/T]\into [(X^\sigma\times \GL_n)/(C_{\GL_n}(\sigma)\times T)]$ induced by the equivariant morphism $x\arr (x,g)$. The diagrams
\[
\begin{tikzcd}
\lbrack X^\sigma/T\rbrack\arrow[rd,swap ,"\pi_{\tilde\sigma}"]\rar["i_{\tilde\sigma}"] & \lbrack(X^\sigma\times \GL_n)/(C_{\GL_n}(\sigma)\times T)\rbrack \dar["\pi '"]\\
& \lbrack X^\sigma /C_{\GL_n}(\sigma)\rbrack 
\end{tikzcd}
\]
are commutative.

We now prove that the composition
\[
\k(X^\sigma\times \GL_n,C_{\GL_n}(\sigma)\times T)_{\tilde\sigma}\xrightarrow{\rho '_*} \k(X,T)_{\sigma '}\xrightarrow{(\rho ')^*} \k(X^\sigma\times \GL_n,C_{\GL_n}(\sigma)\times T)_{\tilde\sigma}
\]
is the multiplication by $\lambda_{-1}((\pi ')^*(\cN))_{\tilde\sigma}$. Since $(i_{\tilde\sigma})_*\colon \k(X^\sigma,T)_{\sigma '}\arr \k(X^\sigma\times \GL_n,C_{\GL_n}(\sigma)\times T)_{\tilde\sigma}$ and $i^*_{\tilde\sigma}\colon \k(X^\sigma\times \GL_n,C_{\GL_n}(\sigma)\times T)_{\tilde\sigma}\arr \k(X^\sigma,T)_{\sigma '}$ are isomorphisms, it is sufficient to prove that for any $x\in \k(X^\sigma,T)_{\sigma '}$ it holds that 
\[
i^*_{\tilde\sigma}((\rho ')^*\rho '_*((i_{\tilde\sigma})_*(x)))=i^*_{\tilde\sigma}((i_{\tilde\sigma})_*(x)\cdot \lambda_{-1}((\pi ')^*(\cN))_{\tilde\sigma}).
\]
Let $\cN_{\sigma '}$ be the conormal bundle of $i_{\tilde\sigma}\colon [X^\sigma/T]\into [(X^\sigma\times \GL_n)/(C_{\GL_n}(\sigma)\times T)] $, and $\cN_T$ be the conormal bundle of $\rho '\circ i_{\tilde\sigma}\colon [X^\sigma/T]\into [X/T]$. Finally, let $\cN '$ be the conormal bundle of $\rho '\colon [(X^\sigma\times \GL_n)/(C_{\GL_n}(\sigma)\times T)]\arr [X/T] $. Since the square of the above diagram is cartesian, we have $(\pi ')^*\cN=\cN '$, so there is an exact sequence
\[
0\arr \cN_{\sigma '}\arr \cN_T\arr i^*_{\tilde\sigma}(\pi ')^*(\cN)\arr 0
\]
which implies that $\lambda_{-1}(i^*_{\tilde\sigma}(\pi ')^*(\cN))\cdot \lambda_{-1}(\cN_{\sigma '})=\lambda_{-1}(\cN_T)$.\\

Using the self-intersection formula twice, we then see that
\[
i^*_{\tilde\sigma}((\rho ')^*\rho '_*((i_{\tilde\sigma})_*(x)))=x\cdot \lambda_{-1}(\cN_T)_{\sigma '}=x\cdot \lambda_{-1}(\cN_{\sigma '})_{\sigma '} \cdot \lambda_{-1}(\pi_{\sigma '}^*(\cN))_{\sigma '}=i^*_{\tilde\sigma}(i_{\tilde\sigma})_*(x)\cdot i^*_{\tilde\sigma}(\pi ')^* (\lambda_{-1}(\cN))_{\sigma }
\]
as we wanted.\\

Now, we note that $W=S_n$ acts on the set of the $\tilde\sigma$'s, and the stabilizer of each of them is $\Delta_\sigma$, the Weyl group of $C_{\GL_n}(\sigma)$. Moreover $W$ acts on the set of the dual cyclic subgroups $\sigma '\subset T$ that are in the same ($\GL_n$ - ) conjugacy class of $\sigma$, and the stabilizer of this action is $\Gamma_\sigma$, the normalizer of $\sigma$ in $\GL_n$. Finally we saw that there is an exact sequence
\[
0\arr \Delta_\sigma\arr\Gamma_\sigma\arr w_{\GL_n}(\sigma)\arr 0
\] 

Now we want to compute the composition of
\[
\Bigl(\underset{\phi\in\text{Aut}(\sigma)/w_{\GL_n}(\sigma)}{\bigoplus}\k(X^{\phi(\sigma)},C_{\GL_n}(\sigma))_\sigma\Bigr)^{\text{Aut}(\mmu_r)}\xrightarrow{\rho_*}\k(X,\GL_n)_\sigma
\]
with
\[
\k(X,\GL_n)_\sigma\xrightarrow{\lambda_{-1}(\cN)^{-1}\cdot\rho^*}
\Bigl(\underset{\phi\in\text{Aut}(\sigma)/w_{\GL_n}(\sigma)}{\bigoplus}\k(X^{\phi(\sigma)},C_{\GL_n}(\sigma))_\sigma\Bigr)^{\text{Aut}(\mmu_r)}.
\]
By the same argument of the case $\cX=[X/T]$, the above composition sends $(x,x,\dots)$ to $\bigl(\frac{\phi(r)}{|w_{\GL_n}(\sigma)|}\cdot \lambda_{-1}(\cN)^{-1}\cdot\rho^*\rho_*(x),\frac{\phi(r)}{|w_{\GL_n}(\sigma)|}\cdot \lambda_{-1}(\cN)^{-1}\cdot\rho^*\rho_*(x),\dots\bigr)$, where $\rho^*\rho_*$ denotes the composition
\[
\k(X^\sigma, C_{\GL_n}(\sigma))^{w_{\GL_n}(\sigma)}_\sigma\xrightarrow{\rho_*}\k(X,\GL_n)_\sigma \xrightarrow{\lambda_{-1}(\cN)^{-1}\cdot\rho^*}\k(X^\sigma, C_{\GL_n}(\sigma))^{w_{\GL_n}(\sigma)}_\sigma. 
\]
To compute this, let us consider the commutative diagrams
\[
\begin{tikzcd}
\k(X^\sigma, C_{\GL_n}(\sigma))^{w_{\GL_n}(\sigma)}_\sigma\rar["\rho_*"]\dar["(\pi ')^*"] & \k(X,\GL_n)_\sigma \dar["\pi^*"]\\
\Bigl( \underset{g\in W/\Delta_\sigma}{\bigoplus}\k(X^{\phi(\sigma)}\times\GL_n,C_{\GL_n}(\sigma)\times T)_{\tilde\sigma}^{\Gamma_\sigma}\Bigr)^W\rar["\rho '_*"]& \Bigl(\underset{g\in W/\Gamma_\sigma}{\bigoplus}\k(X,T)_{\sigma '}\Bigr)^W
\end{tikzcd}
\]
and
\[
\begin{tikzcd}[column sep=15ex]
\k(X,\GL_n)_\sigma \dar["\pi^*"]\rar["\lambda_{-1}(\cN)^{-1}\cdot\rho^*"] & 
 \k(X^\sigma, C_{\GL_n}(\sigma))^{w_{\GL_n}(\sigma)}_\sigma \dar["(\pi ')^*"]\\
\Bigl(\underset{g\in W/\Gamma_\sigma}{\bigoplus}\k(X,T)_{\sigma '}\Bigr)^W\rar[" \lambda_{-1}(\cN)^{-1}\cdot\rho^* "] &  \Bigl(\underset{g\in W/\Delta_\sigma}{\bigoplus}\k(X^{\phi(\sigma)}\times\GL_n,C_{\GL_n}(\sigma)\times T)_{\tilde\sigma}^{\Gamma_\sigma}\Bigr)^W
\end{tikzcd}
\]
Fix any $g\in W/\Delta_\delta$ (recalling $\tilde\sigma=(\sigma,\sigma ')=(\sigma,g^{-1}\sigma g)$) the map $\rho '_*$ sends  
\[
\underset{\phi\in g w_{\GL_n}(\sigma)}{\bigoplus}\k(X^{\phi(\sigma)}\times\GL_n,C_{\GL_n}(\sigma)\times T)_{\tilde\sigma}^{\Gamma_\sigma}\arr \k(X,T)_{\sigma '}
\]
The map $(\pi ')^*\lambda_{-1}(\cN)^{-1}\rho^*\rho_*$ takes an element $(x,\dots,\phi(x),\dots)_{\phi\in g w_{\GL_n}(\sigma)}$; but when $x$ is $\Gamma_\sigma$-invariant it is of the form $(x,x,\dots)$, so the above map is $|w_{\GL_n}(\sigma)|$ times the composition
\[
\k(X^\sigma\times \GL_n,C_{\GL_n}(\sigma)\times T)_{\tilde\sigma}^{\Gamma_\sigma}\xrightarrow{\rho '_*} \k(X,T)^{\Gamma_\sigma}_{\sigma '}\xrightarrow{(\rho ')^*} \k(X^\sigma\times \GL_n,C_{\GL_n}(\sigma)\times T)_{\tilde\sigma}^{\Gamma_\sigma}
\]
We have showed that this composition is the multiplication by $\lambda_{-1}((\pi ')^*(\cN))_{\tilde\sigma}$, thus the map
\[
\k(X^\sigma, C_{\GL_n}(\sigma))^{w_{\GL_n}(\sigma)}_\sigma\xrightarrow{\rho_*}\k(X,\GL_n)_\sigma \xrightarrow{\lambda_{-1}(\cN)^{-1}\cdot\rho^*}\k(X^\sigma, C_{\GL_n}(\sigma))^{w_{\GL_n}(\sigma)}_\sigma 
\]
is a multiplication by $|w_{\GL_n}(\sigma)|$. Given $x\in \k(X^\sigma, C_{\GL_n}(\sigma))^{w_{\GL_n}(\sigma)}_\sigma$, we conclude that
\[
\bigl(\frac{\phi(r)}{|w_{\GL_n}(\sigma)|}\cdot \lambda_{-1}(\cN)^{-1}\cdot\rho^*\rho_*(x),\frac{\phi(r)}{|w_{\GL_n}(\sigma)|}\cdot \lambda_{-1}(\cN)^{-1}\cdot\rho^*\rho_*(x),\dots\bigr)
\]
is equal to $(\phi(r)\cdot x,\phi(r)\cdot x,\dots)$, as we wanted to prove.
\end{proof}

\bigskip

\begin{lemma}
Let $\cX=[X/G]$ and $\sigma\simeq\mmu_r$ be a dual cyclic subgroup; let $\cI_\sigma\cX:=[X^\sigma/C_G(\sigma)]$. Then the composition
\[
\kg{\cI_\sigma\cX}\otimes\QQ(\zeta_r)\into \k(\cI_\sigma\cX)\otimes\R\mmu_r\xrightarrow{\beta_\cX}\k(\cI_\sigma\cX)_\sigma\xrightarrow{\alpha^*_\cX}
\k(\cI_\sigma\cX)\otimes\R\mmu_r\arr \kg{\cI_\sigma\cX}\otimes\QQ(\zeta_r)
\]
is a multiplication by $\frac{1}{r}$. 
\end{lemma}

\begin{proof}
Suppose first that $\cX$ is of the form $[X/G]$, where $X$ is a scheme and $G$ a finite group over $Spec(R)$ containing $\sigma$. Consider the stack $\cY:= [X^\sigma/\sigma]$; then we have an obvious map $\cY\arr \cX$, which induces a map $\cI_\sigma\cY\arr \cI_\sigma\cX$. The latter is the finite covering $[X^\sigma/\sigma]\xrightarrow{\pi} [X^\sigma/C_G(\sigma)]$ by the group $\Gamma:=C_G(\sigma)/\sigma$. By Proposition \ref{prop:descent-geometric-K} the push-forward $\pi_*$ is a surjection in geometric K-theory. 

The stack $\cI_\sigma\cY=X^\sigma\times \cB\sigma$ is the classifying stack of the constant group $\sigma$ over $S:=X^\sigma$. In particular there is a surjection $\pi_*\colon \k(S)\simeq (\R\sigma_S)_\textbf{g}\arr \k(X^\sigma,C_G(\sigma))_\textbf{g}$.
 
We get a commutative diagram
\[
\begin{tikzcd}[sep=small]
\QQ(\zeta_r)\otimes (\R\sigma_S)_\textbf{g}\rar\dar & \R\mmu_r\otimes \R \sigma_S\rar["\beta_{\cB\sigma}"]\dar & (\R \sigma_S)_\sigma \rar["\alpha^*_{\cB\sigma}"]\dar & \R\mmu_r\otimes\R\sigma_S\rar\dar &\QQ(\zeta_r)\otimes (\R\sigma_S)_\textbf{g}\dar \\
\QQ(\zeta_r)\otimes \kg{\cI_\sigma\cX}\rar & \R\mmu_r\otimes \k(\cI_\sigma\cX)\rar["\beta_{\cX}"] & \k(\cI_\sigma\cX)_\sigma \rar["\alpha^*_{\cX}"] & \R\mmu_r\otimes\k(\cI_\sigma\cX)\rar &\QQ(\zeta_r)\otimes \kg{\cI_\sigma\cX}
\end{tikzcd}
\]
To see that the right half commutes, note that the push-forward $\pi_*$ of an element of the form $\cE\otimes V$, where $\cE$ is a sheaf on $X^\sigma$ and $V\in \R\sigma_R$, is just the virtual sheaf $\cE\otimes (\ind_\sigma^{C_G(\sigma)}V)$ seen as a $C_G(\sigma)$-equivariant sheaf on $X^\sigma$; then we may just check the commutativity when $X=R$, where it is clear.

As we have already seen, the composition of the top rows is 
\[
\k(S)\otimes\QQ(\zeta_r)\into\R\mmu_{r,S}\simeq R\mmu_{r,S}\otimes\frac{1+x+\dots+x^{r-1}}{r} \xrightarrow{\Delta}\R\mmu_{r,S}
\]
where $\Delta$ sends $x^i\otimes x^j$ to $\delta_{ij}x^i$, and is equal to $\frac{1}{r}$ times the canonical immersion $\k(S)\otimes\QQ(\zeta_r)\into\R\mmu_{r,S}$, whence our claim holds - trivially - for the top row. But the map $(\R\sigma_S)_\textbf{g}\arr \k(X^\sigma,C_G(\sigma))_\textbf{g}$ is a surjection, so it also holds for the bottom row and the claim is established.\\

Suppose now that $\cX$ is Deligne-Mumford or tame. We use the following lemma (see Appendix C), based on a technique of Toen (\cite{To}, Proposition $4.9$) and a technical result by Vistoli (\cite{Vistoli1}, Lemma $2.7$):
\begin{lemma}
Let $\cX$ be a DM or tame quotient stack. Then there exists a Chow envelope $p\colon\cF\arr \cX$ with $\cF$ a disjoint union of neutral gerbes, that is a finite representable map such that the push-forward $\cI p_*\colon \cI_{\mmu}\cF\arr\cI_{\mmu}\cX$ is surjective.
\end{lemma}

From the previous discussion it follows that we have the result in the case when $\cX=\cB_S G$ for a finite group-scheme $G$ over a scheme $S$; in particular we have it for $\cF$. Now consider the diagram
\[
\begin{tikzcd}
\QQ(\zeta_r)\otimes\kg{\ci[r]\cF}\rar \dar["\cI p_*"]&\QQ(\zeta_r)\otimes\kg{\ci[r]\cF}\dar["\cI p_*"]\\
\QQ(\zeta_r)\otimes\kg{\ci[r]\cX}\rar & \QQ(\zeta_r)\otimes\kg{\ci[r]\cX}
\end{tikzcd}
\]
which is commutative since $p$ is representable. Since $\cI\cX$ is regular we have a projection formula for $f=\cI p$, that is $f_*f^*(\xi)=\xi\cdot f_*(\cO_{\cI\cF})$. But $f$ is finite and surjective, so that $f_*(\cO_{\cI\cF})$ has everywhere nonzero rank in $K_0(\cI\cX)$. By Proposition \ref{prop:decomposition} it is invertible in $\kg{\cI\cX}$, so we conclude that the vertical maps $f=\cI p_*$ are surjective.

Now, the top horizontal map is a multiplication by $\frac{1}{r}$ and $\cI p_*$ is surjective, so the claim also holds for the bottom map.\\ 
\end{proof}
The combination of the two above lemmas gives immediately the proof of the theorem.\\

We conclude this section mentioning that we can compare our morphism with the Toen-Riemann-Roch map (\cite{To}). Let $\Lambda=\QQ(\zeta_\infty)$ be the field generated by the roots of unity; for any $\ZZ$-module $M$ we denote $M_\Lambda$ the $\Lambda-$module $M\otimes_\ZZ\Lambda$. 

Toen defined a map 
\[
\k(\cX)_\Lambda\longrightarrow \Bigl(\kg{\ci[r]\cX}\Bigr)_\Lambda
\]
by composing $$p^*\colon
\k(\cX)_\Lambda\xrightarrow{\lambda_{-1}(\cN)^{-1}\cdot\rho^*} \k(\cI\cX)_\Lambda\xrightarrow{\alpha_\cX^*}\bigoplus_{\substack{r}}(\k(\ci[r]\cX))_\Lambda\otimes\R\mmu_r\arr \bigoplus_{\substack{r}}(\kg{\ci[r]\cX})_\Lambda\otimes\QQ(\zeta_r) 
$$ with the tautological maps $(\kg{\ci[r]\cX})_\Lambda\otimes\QQ(\zeta_r)\arr (\kg{\ci[r]\cX})_\Lambda$ given by $x\otimes\zeta_r^i\arr x\cdot\zeta_r^i$.\\
 
This morphism is clearly the composition of our Lefschetz-Riemann-Roch morphism $\cL$ with the map
 \[
t\colon ((\kg{\cI_{\mmu_r}\cX})_\Lambda\otimes\QQ(\zeta_r))^{\text{Aut}(\mmu_r)}\arr (\kg{\cI_{\mmu_r}\cX})_{\Lambda}
\]
given by $x\otimes\zeta_r^i\arr \frac{\phi(r)}{r}\cdot x\cdot\zeta_r^i$.\\

With the machinery we have developed, we can immediately give a proof of the covariance of Toen's map: since we already know that $\cL$ is covariant, we just need to verify the following

\begin{proposition}\label{cov}
The map $t$ is covariant with respect to proper push-forwards.
\end{proposition}

We conclude that the composition of $t$ with our Riemann-Roch map is a well-defined covariant map with respect to proper push-forwards of relatively tame morphisms, and it coincides with the Toen-Riemann-Roch map. \\

Before giving the proof, we begin with an observation. Consider two positive integers $r|n$ and the projection map $f\colon \cB \mmu_n\arr \cB \mmu_r$.

As we have seen $\cI f_*\colon \k(\cI\cB\mmu_n)\arr \k(\cI\cB\mmu_r)$ induces a map $\cI f_*\colon \kt{\cI\cB\mmu_n}\arr \kt{\cI\cB\mmu_r}$. Let us consider the identity immersion $\sigma=\text{id}\colon \mmu_n\arr\mmu_n$: it induces via $f$ the identity $\sigma '=\text{id}\colon\mmu_r\arr \mmu_r$, so by the above remark we have a map 
\[
 f_*\colon \QQ(\zeta_n)=(\R\mmu_n)_\sigma\arr(\R\mmu_r)_{\sigma '}=\QQ(\zeta_r).
\]

\begin{lemma}\label{trace}
The map $ f_*\colon \QQ(\zeta_n)\arr\QQ(\zeta_r)$ is equal to $\frac{r}{n}\cdot \emph{tr}$, where $\emph{tr}\colon\QQ(\zeta_n)\arr\QQ(\zeta_r)$ is the field-theoretic trace map. 
\end{lemma}

\begin{proof}
Let us make the following observation: the ring $\R\mmu_n$ is a $\R\mmu_r$-module via the pull-back $f^*\colon \R\mmu_r\arr\R\mmu_n$ and - by the projection formula - $f_*$ is a map of $\R\mmu_r$-modules. Moreover $f^*$ gives by localization the immersion $i\colon \QQ(\zeta_r)\into \QQ(\zeta_n)$; 
thus the $\QQ(\zeta_r)$ action on $\QQ(\zeta_n)$ is induced localizing the action of $\R\mmu_r$. We conclude that the push-forward $f_*\colon \QQ(\zeta_n)\arr\QQ(\zeta_r)$ is $\QQ(\zeta_r)$-linear.

We now recall that the bilinear pairing $\text{tr}\colon \QQ(\zeta_n)\times\QQ(\zeta_n)\arr \QQ(\zeta_r)$ given by $(x,y)\arr\text{tr}(x\cdot y)$ is nondegenerate. In particular the $\QQ(\zeta_r)$-linear functional $f_*$ must be of the form $x\arr\text{tr}(k\cdot x)$ for a suitable $k\in\QQ(\zeta_n)$. What we are left to prove is that $k=\frac{r}{n}$.\\

Let us now consider the case $r=1$. Recall that the immersion $\QQ(\zeta_n)\into \R\mmu_n$ sends an element $p(\zeta_n)$ (where $p$ is a polynomial modulo $\phi_n$) to $\psi_n\cdot p\in\R\mmu_n$, where $\psi_n$ is the unique polynomial (modulo $x^n-1$) such that $\psi_n(\zeta_n^i)=1$ if $(n,i)=1$ and is $0$ otherwise.

The push-forward $f_*\colon\R\mmu_n\arr \QQ$ sends a monomial $x^i$ to $1$ if $n|i$ and to $0$ otherwise (this is obviously well-defined modulo $x^n-1$). This is equivalent to the map $p\arr\frac{1}{n}\sum_{i=0}^{n-1}p(\zeta_n^i)$. Combining the definitions, we see that $ f_*\colon \QQ(\zeta_n)\arr\QQ$ sends $p(\zeta_n)$ to 
\[
\frac{1}{n}\sum_{i=0}^{n-1}\psi_n(\zeta_n^i)p(\zeta_n^i)=\frac{1}{n}\sum_{(n,i)=1}p(\zeta_n^i)=\frac{1}{n}\text{tr}(p(\zeta_n)).
\]

Let us return to the general case. We have $f_*\colon \QQ(\zeta_n)\arr\QQ(\zeta_r)=\text{tr}(k\cdot  -)$; in particular the composition $\QQ(\zeta_n)\arr\QQ(\zeta_r)\arr\QQ$ is equal to 

\[
\frac{1}{r}\cdot\text{tr}^{\QQ(\zeta_r)}_\QQ\circ\text{tr}^{\QQ(\zeta_n)}_{\QQ(\zeta_r)}(k\cdot -)=\frac{1}{r}\cdot\text{tr}^{\QQ(\zeta_n)}_\QQ(k\cdot -).
\]
However it is also equal to $\frac{1}{n}\cdot \text{tr}^{\QQ(\zeta_n)}_\QQ$ and by the nondegeneracy of the trace form we conclude that $k=\frac{r}{n}$, as we wanted.
\end{proof}
We can now prove Proposition~\ref{cov}:

\begin{proof}
Let us recall how the push-forward is defined on the left-hand side.
Suppose that we have a proper map $f$ from $\cY=[Y/H]$ to $\cX=[X/G]$ where $X,Y$ are schemes and $G,H$ are general linear groups.  If $Z = X \times_{\cX}Y$, there is an action of $G \times H$ on it, and $[Z/G\times H] = \cY$. The projection $Z \arr X$ is equivariant with respect to $\pi\colon G\times H\arr G$ and gives
   \[
   f\colon\cY = [Z/ G\times H] \arr [X/G] = \cX.
   \]

Let $\rho\simeq\mmu_n\into G\times H$ be a dual cyclic subgroup, whose projections onto $G,H$ are $\sigma_G,\sigma_H$ respectively. Then $C_{G\times H}(\rho)=C_G(\sigma_G)\times C_H(\sigma_H)$.

Let $\sigma\simeq\mmu_r=\sigma_G$; the maps $Z\arr X$ and $\cB_R\rho\arr\cB_R\sigma$ induce a morphism
\[
[Z^\rho/C_G(\sigma_G)\times C_H(\sigma_H)]\times \cB_R\rho\arr [X^\sigma/C_G(\sigma)]\times\cB_R\sigma
\]
and this is the restriction of $\widetilde{\ci{f}}$ to the component of $\widetilde{\ci\cY}$ relative to $\rho$. As we have seen it induces a map $(\kg{\cI_{\mmu_n}\cY})_\Lambda\otimes\QQ(\zeta_n)\arr (\kg{\cI_{\mmu_r}\cX})_\Lambda\otimes\QQ(\zeta_r)$.

Let $\Gamma=\text{Gal}(\QQ(\zeta_n):\QQ(\zeta_r))$. Then acting on $\rho$ by elements of $\Gamma$ gives all the conjugacy classes of dual cyclic subgroups which lie over $\sigma$. Let $\Gamma_1$ be the stabilizer of $\rho$ in $\Gamma$ (that is the subgroup of elements $\gamma$ such that $\gamma(\rho)$ and $\rho$ are conjugated). Let $\cI_\rho\cY:=[Z^\rho/C_G(\sigma_G)\times C_H(\sigma_H)]$ and $\cI_\sigma\cX:=[X^\sigma/C_G(\sigma)]$; we want to prove that the following diagram commutes:
\[
\begin{tikzcd}
\Bigl(\underset{\gamma\in\Gamma/\Gamma_1}{\bigoplus}(\kg{\cI_{\gamma(\rho)}\cY})_\Lambda\otimes\QQ(\zeta_n)\Bigr)^{\text{Aut}(\mmu_n)}\rar\dar &
\underset{\gamma\in\Gamma/\Gamma_1}{\bigoplus}(\kg{\cI_{\gamma(\rho)}\cY})_\Lambda\dar \\
(\kg{\cI_{\sigma}\cX})_\Lambda\otimes\QQ(\zeta_r)\rar &(\kg{\cI_{\sigma}\cX})_\Lambda
\end{tikzcd}
\]
Any element of the top left box is of the form $\tilde x=(x,\dots,\gamma(x),\dots)$ for $x\in ((\kg{\cI_\rho\cY})_\Lambda\otimes\QQ(\zeta_n))^{\Gamma_1}$.

Let $x=\sum x_i\otimes\zeta_n^i$. All the elements $(\gamma(x))$ have the same image in $(\kg{\cI_{\sigma}\cX})_\Lambda\otimes\QQ(\zeta_r)$, since $\gamma$ has trivial image in $\text{Aut}(\mmu_r)$. Thus the composition 
\[
\Bigl(\underset{\gamma\in\Gamma/\Gamma_1}{\bigoplus}(\kg{\cI_{\gamma(\rho)}\cY})_\Lambda\otimes\QQ(\zeta_n)\Bigr)^{\text{Aut}(\mmu_n)}\arr
(\kg{\cI_{\sigma}\cX})_\Lambda\otimes\QQ(\zeta_r)\arr(\kg{\cI_{\sigma}\cX})_\Lambda
\]
sends $\tilde x$ to 
\[
|\Gamma/\Gamma_1|\cdot\sum_i\bigl(\cI f_*(x_i)\cdot\frac{\phi(r)}{r}\cdot\frac{r}{n}\text{tr}(\zeta_n^i)\bigr)
=|\Gamma/\Gamma_1|\cdot\frac{\phi(r)}{n}\cdot\sum_i\sum_{\gamma\in\Gamma}\cI f_*(x_i)\cdot \gamma(\zeta_n^i)
\]
However, as we observed before, $\cI f_*(x_i)=\gamma(\cI f_*(x_i))=\ci f_*(\gamma(x_i))$ so the above term is equal to
\[
\begin{aligned}
|\Gamma/\Gamma_1|\cdot\frac{\phi(r)}{n}\cdot\sum_i\sum_{\gamma\in\Gamma}\cI f_*(\gamma(x_i))\cdot \gamma(\zeta_n^i)&=&|\Gamma/\Gamma_1|\cdot\frac{\phi(r)}{n}\cdot\sum_{\gamma\in\Gamma/\Gamma_1}\sum_{\gamma_1\in\Gamma_1}\gamma\Bigl(\sum_i(\cI f_*(\gamma_1(x_i\cdot\zeta_n^i)))\Bigr) \\
{}&=&
|\Gamma/\Gamma_1|\cdot\frac{\phi(r)}{n}\cdot|\Gamma_1|\cdot\sum_{\gamma\in\Gamma/\Gamma_1}\gamma\Bigl(\sum_i\Bigl(\cI f_*(x_i)\cdot \zeta_n^i\Bigr)\Bigr)
\end{aligned}
\]
recalling that $\gamma_1(x)=x$ for any $\gamma_1\in\Gamma_1$. But $|\Gamma|=\frac{\phi(n)}{\phi(r)}$ so we can rewrite this as
\[
\frac{\phi(n)}{n}\cdot\sum_{\gamma\in\Gamma/\Gamma_1}\gamma\Bigl(\sum_i\Bigl(\cI f_*(x_i)\cdot\zeta_n^i\Bigr)\Bigr)
\]
that is the image of $\tilde x$ under the composition
\[
\Bigl(\underset{\gamma\in\Gamma/\Gamma_1}{\bigoplus}(\kg{\cI_{\gamma(\rho)}\cY})_\Lambda\otimes\QQ(\zeta_n)\Bigr)^{\text{Aut}(\mmu_n)}\arr
\underset{\gamma\in\Gamma/\Gamma_1}{\bigoplus}(\kg{\cI_{\gamma(\rho)}\cY})_\Lambda\arr (\kg{\cI_{\sigma}\cX})_\Lambda.
\]
The Proposition is thus proved.
\end{proof}

\section{The rational Toen-Riemann-Roch map}

Let $I$ be the moduli space of $\cI\cX$. Using the results of the previous section, we can produce a covariant isomorphism
\[
\k(\cX)\longrightarrow A^*(I).
\]
with image in the Chow groups of $I$. This map, with respect to the classical Toen-Riemann-Roch morphism, has the advantage of having $\QQ$-coefficients and the disadvantage of being non-canonical.\\

We start by providing a $\QQ$-isomorphism
\[
\k(\cX)\longrightarrow \kg{\cI\cX}.
\]

In general, given a finite group $G$ and a finitely generated $k[G]-$module $M$, there is a canonical isomorphism 
\[
M\longrightarrow (M\otimes k[G])^G
\]
given by $m\arr \frac{1}{|G|}\underset{g\in G}{\sum}g(m)\otimes g$.\\

In particular for any $r$ we have an isomorphism 
\[
\kg{\ci[r]\cX}\arr \Bigl(\kg{\ci[r]\cX}\otimes\QQ(\zeta_r)\Bigr)^{\text{Aut}(\mmu_r)}
\]
once we can provide an isomorphism $\QQ(\zeta_r)\simeq k[\text{Aut}(\mmu_r)]$. This is given by a choice of a \emph{normal basis} $(x,\dots, \gamma(x),\dots )_{\gamma\in\text{Aut}(\mmu_r)}$ for $\QQ(\zeta_r)$.\\

These maps combine into an isomorphism $\k(\cX)\longrightarrow \kg{\cI\cX}$ for any $\cX$. In order to make these maps covariant for proper push-forwards, we invoke the following result of Lenstra (\cite{Len}):

\begin{proposition}
There exist, for any $n$, an element $x_n\in\QQ(\zeta_n)$ such that
\begin{enumerate}
\item The set $\{\gamma(x_n)\}_{\gamma\in\text{Aut}(\mmu_r)}$ is a normal basis for $\QQ(\zeta_n)$.
\item For any $n,m$ with $m|n$ we have $\text{Tr}^{\QQ(\zeta_n)}_{\QQ(\zeta_m)}(x_n)=x_m$.

\end{enumerate}
\end{proposition}

Having this result, we can now take the normal bases $\{\gamma\bigl(r\cdot x_r\bigr)\}_{\gamma\in\text{Aut}(\mmu_r)}$ in our construction of the map $$\phi\colon \kg{\ci[r]\cX}\arr \Bigl(\kg{\ci[r]\cX}\otimes\QQ(\zeta_r)\Bigr)^{\text{Aut}(\mmu_r)}.$$

\begin{proposition}
The map $\phi$ is covariant with respect to proper push-forward.
\end{proposition}

\begin{proof}
The proof is very similar to the proof of Proposition~\ref{cov}. Suppose that we have a proper map $f$ from $\cY=[Y/H]$ to $\cX=[X/G]$.  As usual, if $Z = X \times_{\cX}Y$, there is an action of $G \times H$ on it, and $[Z/G\times H] = \cY$. The projection $Z \arr X$ is equivariant with respect to $\pi\colon G\times H\arr G$ and gives
   \[
   f\colon\cY = [Z/ G\times H] \arr [X/G] = \cX.
   \]

Let $\rho\simeq\mmu_n\into G\times H$ be a dual cyclic subgroup, whose projections onto $G,H$ are $\sigma_G,\sigma_H$ respectively. Then $C_{G\times H}(\rho)=C_G(\sigma_G)\times C_H(\sigma_H)$.

Let $\sigma\simeq\mmu_r=\sigma_G$ and suppose that $x\in\kg{{\ci}_\rho\cY}$, so that $\cI f_*(x)\in \kg{{\ci}_\sigma\cX}$. Let $n=|\rho|$ and $r=|\sigma|$.

Then we want to prove that the diagram
\[
\begin{tikzcd}
\kg{\ci[n]\cY}\rar["\phi"]\dar["\cI f_*"]& \kg{\ci[n]\cY}\otimes\QQ(\zeta_n)\dar["\tilde{\cI f}_*"]\\
\kg{\ci[r]\cX}\rar["\phi"] &\kg{\ci[r]\cX}\otimes\QQ(\zeta_r)
\end{tikzcd}
\]
is commutative. The composition of the bottom horizontal and left vertical arrows send $$x\arr \frac{r}{\phi(r)}\underset{\gamma\in \text{Aut}(\mmu_r)}{\sum}\gamma(\cI f_*(x))\otimes\gamma(x_r).$$

On the other hand, using Lemma~\ref{trace} we immediately see that the composition of the other arrows send $$ x \arr \frac{n}{\phi(n)}\underset{\gamma\in \text{Aut}(\mmu_n)}{\sum}\cI f_*(\gamma(x))\otimes\frac{r}{n}\text{Tr}^{\QQ(\zeta_n)}_{\QQ(\zeta_r)}(x_n).$$
Let $\Gamma$ be the Galois group of $\mmu_n$ over $\mmu_r$; then for any $\gamma\in\Gamma$ we have $\cI f_*=\cI f_*\circ\gamma$, so we can rewrite the last expression as 
\[
\frac{n}{\phi(n)}\cdot |\Gamma|\underset{\gamma\in \text{Aut}(\mmu_n)/\Gamma}{\sum}\cI f_*(\gamma(x))\otimes\frac{r}{n}\text{Tr}^{\QQ(\zeta_n)}_{\QQ(\zeta_r)}(x_n)
\]
But $\text{Aut}(\mmu_n)/\Gamma=\text{Aut}(\mmu_r)$ and $\text{Tr}^{\QQ(\zeta_n)}_{\QQ(\zeta_r)}(x_n)=x_r$, so this is equal to $$\frac{r}{\phi(r)}\underset{\gamma\in \text{Aut}(\mmu_r)}{\sum}\gamma(\cI f_*(x))\otimes\gamma(x_r)$$ as we wanted.
\end{proof}

Now, let $I$ be the coarse moduli space of $\cI\cX$. We have an isomorphism $\pi_*\colon \kg{\cI\cX}\arr \kg{I}\simeq \k(I)$, and the usual covariant Riemann-Roch map $\tau_I\colon \k(I)\arr A^*(I)$, where $A^*(-)$ is the Chow group functor. Summarizing, the composition

\[
\k(\cX)\rightarrow \kg{\cI\cX}\xrightarrow{\pi_*} \kg{I}\simeq \k(I)\xrightarrow{\tau_I} A^*(I)
\]
provides a Riemann-Roch map which is covariant and has $\QQ-$coefficients.

\clearpage
\begin{appendices}
\counterwithin{theorem}{section}
\section{Homomorphism schemes from diagonalizable groups}

Let us analyze the structure of $\rh{\Delta}G$ when $\GL_{n_{1}}  \times \dots \times \GL_{n_{r}}$ is a product of general linear groups.

\begin{proposition}\label{prop:cyclotomic-in-GL}
Assume that $G$ is a product of general linear groups. Then $\rh{\Delta}G$ is represented by an open subscheme in a disjoint union of product of Grassmannians. The action of $G$ on each connected component of $\rh{\Delta}G$ is transitive, and each connected component contains a $k$-rational point. Furthermore, $\rhin{\Delta}{G}$ is a union of connected components of $\rh{\Delta}G$.

Given two homomorphism $f$, $g\colon \Delta \arr G$, thought of as $k$-rational points of $\rh{\Delta}G$, then $f$ and $g$ are in the same connected component of $\rh{\Delta}G$ if and only if then are conjugate by an element of $G(k)$.
\end{proposition}

\begin{proof}
Let $G=\GL_{n_1}\times\dots\times\GL_{n_r}$. Clearly $\rhin{\Delta}{G}\subseteq \rh{\Delta}{G}$ is $G-$invariant, so it is a union of connected components of $\rh{\Delta}{G}$ if we know that they are $G-$orbits. We are then left to prove the remaining claims.

A homomorphism $\Delta_S\arr G_S$ corresponds to eigenspaces decompositions $\cO^{n_i}_{S} = \bigoplus_{\chi \in \widehat{\Delta}} V_{\chi,i}$ for every $i=1,\dots, r$. For every function $d\colon \widehat{\Delta}\arr \NN^r$, let $\rH^{d}_{\Delta}(\GL_{n}) \subseteq \rh{\Delta}{\GL_{n}}$ be the subfunctor of the homomorphisms $\Delta\arr G$ such that for all $i$'s the eigenspaces $V_{\chi,i}$ have constant rank $\pi_id(\chi)$ (where $\pi_i$ denotes the $i$th component of $d$). We have a decomposition $\rh{\Delta}{G} = \coprod_{d}\rH_{\Delta}^{d}(G)$. If $0 \leq m \leq n$, denote by $\GG(m,n)$ the Grassmannian of $m$-dimensional subspaces of $k^{n}$. It has a tautological action of $\GL_n$, and there is a $G-$equivariant open embedding
   \[
   \rH_{\Delta}^{d}(\GL_{n}) \subseteq \prod_{\chi,i}\GG\bigl(\pi_id(\chi), n_i\bigr)
   \] 
which identifies the former with the open subscheme of $\prod_{\chi,i}\GG\bigl(\pi_id(\chi), n\bigr)$ such that for every fixed $i$ the subspaces given by $\GG\bigl(\pi_id(\chi), n_i\bigr)$ are in direct sum. The latter is an integral scheme and the $G-$action on it is clearly transitive, proving the first claim.

The spaces $\cO^{n_i}_R$ obviously have a $d$-decomposition for every $d$ such that $\sum_{\chi} \pi_id(\chi)=n_i$, so every connected component has a rational point, and we conclude exactly as above that the $G(R)$ action is transitive on the rational points.
\end{proof}

If $\phi\colon \Delta \arr \Delta'$ is a homomorphism of diagonalizable group schemes, composing with $\phi$ gives a natural transformation $\rh{\Delta'}G \arr \rh{\Delta}G$. Call $\rQ(\Delta)$ the set of quotients of $\Delta$. For each $\Delta' \in \rQ(\Delta)$, consider the composite $\rhin{\Delta'}G \subseteq \rh{\Delta'}G \arr \rh{\Delta}G$, which is immediately seen to be a monomorphism. This induces a $G$-equivariant morphism of Zariski sheaves
   \[
   \coprod_{\Delta' \in \rQ(\Delta)} \rhin{\Delta'}G \arr  \rh{\Delta}G\,.
   \]

\begin{proposition}\call{prop:cyclotomic-in-G}\hfil
\begin{enumerate1}

\itemref{1} The functors $\rh{\Delta}{G}$ and\/ $\rhin{\Delta}G$ are represented by a quasi-projective scheme over $R$.

\itemref{2} If $\Delta' \in \rQ(\Delta)$, then $\rhin{\Delta'}G$ is open and closed in $\rh{\Delta}G$.

\itemref{3} The morphism $\coprod_{\Delta' \in \rQ(\Delta)} \rhin{\Delta'}G \arr  \rh{\Delta}G$ is an isomorphism.

\itemref{4} If $G$ is finite and linearly reductive, then $\rh{\Delta}{G}$ and $\rhin{\Delta}G$ are finite over $\spec R$. If $G$ is étale and $R$ equicharacteristic,  $\rh{\Delta}{G}$ and $\rhin{\Delta}G$ are finite over $\spec R$.


\end{enumerate1}

\end{proposition}

\begin{proof}

Choose an embedding $G \subseteq \GL_{n}$ for some $n$; this gives an embedding of functors $\rh{\Delta}{G} \subseteq \rh{\Delta}{\GL_{n}}$.

\begin{lemma}
The inclusion $\rh{\Delta}{G} \subseteq \rh{\Delta}{\GL_{n}}$ is a closed embedding.
\end{lemma}

\begin{proof}
This is standard.
\end{proof}

Clearly we have $\rhin{\Delta}G = \rh{\Delta}G \cap \rhin{\Delta}{\GL_{n}}$. More generally, if $\Delta' \in \rQ(\Delta)$, the inverse image of $\rh {\Delta}G \subseteq\rhin{\Delta}{\GL_{n}}$ in $\rhin{\Delta'}\GL$ equals $\rhin{\Delta'}G$. Hence, to prove \refpart{prop:cyclotomic-in-G}{1}, \refpart{prop:cyclotomic-in-G}{2} and \refpart{prop:cyclotomic-in-G}{3} we can assume that $G = \GL_{n}$.

So it is enough to prove that $\rh{\Delta}{\GL_{n}}$ is represented by a quasi-projective scheme over $k$. Let $\widehat{\Delta}$ be the group of characters $\Delta \arr \gm$ of $\Delta$. By the standard description of representations of diagonalizable groups, a representation $\Delta_{S} \arr \GL_{n,S}$ corresponds to a decomposition $\cO^{n}_{S} = \bigoplus_{\chi \in \widehat{\Delta}} V_{\chi}$ into eigenspaces. If $d\colon \widehat{\Delta} \arr \NN$ is a function, denote by $\rH^{d}_{\Delta}(\GL_{n}) \subseteq \rh{\Delta}{\GL_{n}}$ the subfunctor of those representations of $\Delta$ such that the corresponding eigenspace $V_{\chi}$ has constant rank $d(\chi)$. We have a decomposition of Zariski sheaves $\rh{\Delta}{\GL_{n}} = \coprod_{d}\rH_{\Delta}^{d}(\GL_{n})$. If $0 \leq m \leq n$, denote by $\GG(m,n)$ the Grassmannian of $m$-dimensional subspaces of $A^{n}$. There is an obvious embedding of functors
   \[
   \rH_{\Delta}^{d}(\GL_{n}) \subseteq \prod_{d}\GG\bigl(d(\chi), n\bigr)
   \] 
which is easily seen to be an open embedding. This proves \refpart{prop:cyclotomic-in-G}{1}.

Furthermore, if $d\colon \widehat{\Delta} \arr \NN$ is a function, denote by $\Delta'_{d}$ the quotient of $\Delta$ such that $\widehat{\Delta'_{d}} \subseteq \widehat{\Delta}$ is the group generated by the $\chi \in \widehat{\Delta}$ with $d(\chi) > 0$. Then it is easy to see that $\rhin{\Delta'}{\GL_{n}} \subseteq \rh{\Delta}{\GL_{n}}$ is the union of the components $\rH_{\Delta}^{d}(\GL_{n})$ with $\Delta'_{d} = \Delta'$. This proves \refpart{prop:cyclotomic-in-G}{2} and \refpart{prop:cyclotomic-in-G}{3}.

To prove \refpart{prop:cyclotomic-in-G}{4}, assume that $G$ is finite and linearly reductive or étale.

If $\Delta = \Delta' \times \Delta''$ is a decomposition into the product of two diagonalizable subgroups, and assume that $\rH_{\Delta'}(G)$ and $\rH_{\Delta''}(G)$ are finite over $R$; let us show that $\rH_{\Delta}(G)$ is also finite. We get an obvious morphism $\rH_{\Delta}(G) \arr \rH_{\Delta'}(G) \times \rH_{\Delta''}(G)$; let us show that this is a closed embedding. In fact, let $S \arr \rH_{\Delta'}(G) \times \rH_{\Delta''}(G)$ be a morphism, corresponding to an object $(f', f'')$ of $\bigl(\rH_{\Delta'}(G) \times \rH_{\Delta''}(G)\bigr)(S)$; here, $f'\colon \Delta'_{S} \arr G_{S}$ and $f''\colon \Delta''_{S} \arr G_{S}$ are homomorphisms of group schemes. Then $(f', f'')$ comes from a (unique) object of $\rH_{\Delta}(G)$ if and only if $f'$ and $f''$ commute, that is, the morphism $(\Delta'\times \Delta'')_{S} \arr G_{S}$ that sends a pair $(\delta', \delta'')$ into $f'(\delta')f''(\delta'')f'(\delta')^{-1}f''(\delta'')^{-1}$ factors through the identity section $S \arr G_{S}$. Then the result follows from the following standard fact.

\begin{lemma}
Let $X \arr S$ and $Y \arr S$ be morphisms of schemes, $f$, $g\colon X \arr Y$ morphisms of $S$-schemes. Assume that $X \arr S$ is finitely presented, finite and flat, while $Y \arr S$ is separated. Then the functor from schemes to sets, sending a scheme $T$ into the set of morphisms $T \arr S$ such that the pullbacks $f_{T}$, $g_{T}\colon  X_{T} \arr Y_{T}$ coincide is represented by a closed subscheme of $S$.
\end{lemma}

Now let us prove the result in general. We may assume that $R$ is a local ring, with residue field $k$. If $G$ is linearly reductive then after extending the base we may assume that $G$ is well-split, that is, a semidirect product $G_{1} \ltimes G_{0}$, where $G_{1}$ is constant, of order not divisible by $\cha k$, and $G_{0}$ is diagonalizable. We can split $\Delta$ into a finite product of group schemes of type $\mmu_{p^{r}}$, where $p$ is a prime; because of the previous step, we can assume that $\Delta = \mmu_{p^{r}}$. 

If $p = \cha k$, suppose first that $G$ is well-split: then $\underhom_{A}(\mmu_{p^{r}}, G_{1}) = \spec A$, so $\rH_{\mmu_{p^{r}}}(G) = \rH_{\mmu_{p^{r}}}(G_{0})$; and, because of Cartier duality, $\rH_{\mmu_{p^{r}}}(G_{0})$ is a finite union of copies of $\spec A$. 
On the other hand if $G$ is étale then $\underhom_{A}(\mmu_{p^{r}}, G) = \spec A$, since $\mmu_{p^{r}}$ is infinitesimal when $A$ is equicharacteristic.

If $p \neq \cha k$, then $\mmu_{p^{r}}$ is a constant cyclic group scheme of order $p^{r}$; hence $\rH_{\mmu_{p^{r}}}(G)$ is represented by the inverse image of the identity $\spec A \subseteq G$ via the map $G \arr G$ defined by $x \arrto x^{p^{r}}$. This ends the proof of \refpart{prop:cyclotomic-in-G}{4}.
\end{proof}

\section{The Weyl group action on torus-equivariant K-theory}

In this section we briefly discuss the isomorphism $\k(X,GL_n)\simeq \k(X,T)^W$, where $T$ is a maximal torus and $W$ the associated Weyl group. Let $G:=GL_n$.\\

To see this, we first note that $\k(X,T)\simeq \k(X,B)$ where $B$ is the associated Borel subgroup of $G$ (since $[X/T]\arr [X/B]$ is an affine bundle).

Now let $E\arr \cX:=[X/G]$ be the associated vector bundle; we recall here how to resume the flag bundle $[X/B]\arr \cX$ as a composition of projective bundles $\PP_i(E)\arr \cX$.

Let $\PP_1(E):=\PP(E)\xrightarrow{\pi_1}\cX$ be the projective bundle associated to $E$. Then we define $E_2$ as the tatutological vector bundle on $\PP_1(E)$ associated to $\PP_1(E)$ (it is a codimension $1$ subbundle of $\pi_1^*E$), and $\PP_2(E):=\PP(E_2)\xrightarrow{\pi_2}\PP_1(E)$; inductively it is then clear how, given $\PP_i(E)\xrightarrow{\pi_{i}} \PP_{i-1}(E)$, we can define $E_{i+1}$ as the tautological bundle of $\PP_i(E)$ (a codimension $1$ subbundle of $\pi_i^*E_i$) and $\PP_{i+1}(E):=\PP(E_{i+1})$. At the end of this process we get $\PP_n(E)=[X/B]$.\\

By the projective bundle theorem in equivariant K-theory (see \cite{Thomason2}), we have isomorphisms $\k(\PP_{i+1}(E))\simeq\prod_{k=0}^{n-i}\k(\PP_i(E))$, and by induction we get that $\k([X/B])$ is isomorphic to a product of $n!$ factors, each isomorphic to $\k(\cX)$; The pull-back identifies $\k(\cX)$ as the diagonal subgroup.

These factors correspond to the pieces of the Bruhat decomposition of $G/B$, and the induced action of the Weyl group is just the tautological permutation of the factors. Then it is clear that $\k(X,G)\simeq \k(X,B)^W\simeq \k(X,T)^W$.

\section{Existence of envelopes for DM and tame stacks}

In this appendix we want to prove the existence of finite envelopes by unions of neutral gerbes for reduced DM and tame algebraic stacks.\\

Let $\cX$ be such a stack, over a fixed base ring $R$, and let $M$ be its space of moduli. By \cite{Vistoli2}, Theorem 2.7, there exists a finite surjective cover $U\arr\cX$ from a scheme. Let us consider the cartesian square
\[
\begin{tikzcd}
\cF\ar[r]\ar[d] & U\ar[d]\\
\cX\ar[r] & M
\end{tikzcd}
\]
where $U\arr M$ is the composition $U\arr \cX\arr M$. It is immediate, by \cite{dan-olsson-vistoli}, Corollary 3.3, to see that $U$ is the moduli space of $\cF$.

Let us take the normalization $\overline{\cF}$ of $\cF$ (which exists by \cite{nor}), and its moduli space $\overline{U}$, that is again normal (since a categorical quotient of a normal variety is normal). We want to prove that $\overline{\cF}$ is a neutral gerbe over $\overline{U}$; this is sufficient to conclude. Indeed, the map $\overline{\cF}\arr \cX$ admits a fppf-local section over an open dense substack $\cV$ of $\cX$ that is a gerbe. In particular, the map $\cI_{\mmu}\cF\arr \cI_{\mmu}\cX$ is generically surjective. By noetherian induction we can ensure that the closed substack $\cX ':= \cX\backslash \cV$ admits a finite envelope $\cF '$, with $\cF'$ a finite union of neutral gerbes; then $\overline{\cF} \sqcup \cF '$ is the sought envelope.\\

Thus we just need to prove that $\overline{\cF}\arr \overline{U}$ is a gerbe. 

Consider the commuting square
\[
\begin{tikzcd}
\overline{\cF}\ar[r]\ar[d] & \overline{U}\ar[d]\\
\cF \ar[r] & U \ar[l,dashed,bend left=25]
\end{tikzcd}
\]
and note that by construction the arrow $\cF\arr U$ admits a section. From the universal property of normalization, the top row $\overline{\cF}\arr \overline{U}$ also admits a section. More precisely, take a smooth atlas $C\arr \cF$ and define $V:=C\times_\cF \overline U$. Consider the following diagram with cartesian quadrangles:
\[
\begin{tikzcd}
\overline{C\times_\cF C} \ar[d]\ar[r,shift left]\ar[r,shift right] &\overline{C} \ar[r]\ar[d] & \overline{\cF}\ar[d]\\
C\times_\cF C \ar[r,shift left]\ar[r,shift right] &C \ar[r] & \cF \\
V\times_{\overline U}V \ar[r,shift left]\ar[r,shift right]\ar[u] &V \ar[u]\ar[r] & \overline U \ar[u]
\end{tikzcd}
\]
By definition the top horizontal row is a presentation for $\overline{\cF}$ and the two bottom-left vertical arrows are dominant, hence by the universal property of normalization all the bottom vertical arrows admit a lift to the top row. In particular the moduli space map $\overline{\cF}\arr \overline U$ admits a section $\overline U\arr \overline{\cF}$.\\

 We conclude with the following

\begin{lemma}
Let $\cF$ be a normal DM or tame algebraic stack such that the projection $\cF\arr U$ to its moduli space admits a section. Then $\cF$ is a neutral gerbe.
\end{lemma}

\begin{proof}
First of all, we note that we can take a rigidification $\cF\arr \cG$; in other words we can exhibit $\cF$ as a gerbe over a generically trivial stack $\cG$ (see the Appendix to \cite{dan-olsson-vistoli} for a reference). Actually, to prove the Lemma we can work étale-locally on $U$, and thus we may suppose that $\cF$ is connected and of the form $[X/G]$, where $G$ is a finite étale or linearly reductive group (see \cite{LM}, Theorem $6.1$ or \cite{dan-olsson-vistoli}, Theorem $3.2$, respectively). In that case we can explicitly construct the rigidification.

Passing to a further étale cover, we may take $G=P\rtimes T$, where $T$ is discrete and $P$ of multiplicative type (\cite{dan-olsson-vistoli}, Lemma $2.17$). Moreover we can make some simplifying assumptions:

\begin{enumerate}
\item  First, we can assume that $G$ is either discrete or of multiplicative type. Indeed, let $\cF=[X/G]$ and $\cF_{\textbf{dm}}=[(X/P)/(G/P)]$ (where the square bracket indicates the stacky quotient and the round bracket the schematic quotient). Consider the cartesian diagram
\[
\begin{tikzcd}
\cU \ar[r]\ar[d] & U \ar[d]\ar[ld] \\
\cF \ar[r] & \cF_{\textbf{dm}}
\end{tikzcd}
\]
where the diagonal arrow is the section $U\arr \cF$. Then $\cF_{\textbf{dm}}$ has only étale stabilizers; moreover $U$ is the moduli space of $\cU$ and $\cF_{\textbf{dm}}$ and by construction in both cases the projection to $U$ admits a section. In particular if we know the result for stacks that have either étale or diagonalizable stabilizers we can deduce the general case: indeed we can infer  that $\cF_{\textbf{dm}}\arr U$ is a gerbe and $\cF\arr \cF_{\textbf{dm}}$ is also a gerbe; for the second, it is sufficient to show that its base-change $\cU\arr U$ is a gerbe, which we infer looking at the cartesian cube
\[
\begin{tikzcd}
& \cV \ar[ld] \ar[rr] \ar[dd] & &  V \ar[dd] \ar[ld]\\
\cU \ar[rr]\ar[dd] & & U \ar[dd] &  \\
& \lbrack X\slash P\rbrack \ar[ld]\ar[rr] & & (X/P) \ar[ld] \\
\cF \ar[rr] & & \cF_{\textbf{dm}} &
\end{tikzcd}
\]
Indeed $\cV$, being a base-change of $[X/P]$, has only diagonalizable stabilizers and $\cV\arr V$ by construction admits a section, so we can deduce that it is a gerbe; this implies that $\cU\arr U$ is also a gerbe.

If  $\cF_{\textbf{dm}}\arr U$ and $\cF\arr \cF_{\textbf{dm}}$ are both gerbes, then $\cF\arr U$ is also a gerbe, which must then be trivial.
\item  Second, we note that we may assume that $X$ is connected. Indeed, let $X$ be a disjoint union of connected components $X=X_1\sqcup\dots\sqcup X_n$; since $[X/G]$ is connected, the map $X_1\times G\xrightarrow{m} X$ is surjective. Let $G_1\leq G$ be the maximal subgroup sending $X_1$ to itself, defined by the cartesian square
\[
\begin{tikzcd}
Z=X_1\times G_1 \ar[r]\ar[d] & X_1\ar[d]\\
X_1\times G \ar[r,"m"] & X
\end{tikzcd}
\]
More precisely, it is readily checked that $Z$ is an open and closed subgroup functor of $X_1\times G$; moreover a flat subgroup of a group which is constant or of multiplicative type is also of the same type; if the base is connected then it must come from a subgroup of $G$, which we call $G_1$.

Then $X_1\times G\xrightarrow{m} X$ is a $G_1$-principal bundle, so that $X_1\times_{G_1}G\arr X$ is an isomorphism. Taking $G-$quotients we get that 
\[
[X_1/G_1]\arr [X/G]
\]
is an isomorphism, so me can assume without loss of generality that $X=X_1$.
\end{enumerate}

Now let $H$ be the stabilizer of a generic point of $X$; we may assume that $H$ is a normal subgroup of $G$. Indeed, if $G$ is of multiplicative type this follows from the fact that the stabilizer of any generic point of $X$ must be of multiplicative type and come from a subgroup of $G$; moreover $G$ is abelian, so $H$ is certainly normal. 

If $G$ is discrete then $X$ is integral (as it is étale over a normal stack, hence normal itself), so $H$ is the stabilizer of the unique generic point of $X$ (and as before comes from a subgroup of $G$); again, $H$ must be normal in $G$.

Let now $Z\subseteq X$ be the fixed locus of $H$, which is closed. Then we have that $Z$ is $G$-invariant, since $H$ is normal in $G$; but $[X/G]$ is reduced (since it is normal) and $[Z/G]$ is a closed substack containing the generic point of $|X/G|$: this implies that $[Z/G]=[X/G]$. We conclude that $Z=X$, as we wanted.\\ 

Summarizing, $H$ acts trivially on $X$, so there is a morphism
\[
[X/G]\arr [X/(G/H)]
\] 
This is easily seen to be the rigidification map.\\

Still working étale-locally and maintaining the above notation we see that the rigidification $\cG$ is again normal. Indeed, assuming again without loss of generality that $H$ is normal, we have $\cG=[X/(G/H)]$ when $\cF=[X/G]$; taking a faithful representation $G\arr \GL_n$ we can write $\cG=[(X\times_G\GL_n)/(\GL_n/H)]$. But $\GL_n/H$ is smooth and thus the map $X\times_G\GL_n\arr \cG$ is smooth; the former scheme is normal by assumption and we conclude that $\cG$ is normal as well.

Then we have that $\cG\arr U$ is generically an isomorphism and admits a section, that must be finite and representable. Since $\cG$ is normal we conclude that $U\arr \cG$ is an isomorphism. In particular $\cG$ is a scheme.

Finally, $\cF\arr U$ is a gerbe over a scheme with a section, so it is a neutral gerbe, as we wanted.
\end{proof}

\end{appendices}

\end{document}